\documentclass[a4paper,reqno]{amsart}

\usepackage{amsfonts,amssymb,amscd,amsmath,latexsym,amsbsy,enumitem,stmaryrd,a4wide,verbatim,color,subfiles,placeins,mathtools}
\usepackage{mwe}
\usepackage{todonotes}

\usepackage{hyperref}
\usepackage{tikz}
\usetikzlibrary{arrows,positioning,shapes.misc}
\theoremstyle{plain}
\newtheorem{theorem}{Theorem}[section]
\newtheorem{corollary}[theorem]{Corollary}

\newtheorem{lemma}[theorem]{Lemma}
\newtheorem{proposition}[theorem]{Proposition}

\theoremstyle{remark}
\newtheorem{remark}[theorem]{Remark}
\numberwithin{equation}{section}

\makeatletter
\def\subsection{\@startsection{subsection}{2}%
  \z@{.7\linespacing\@plus\linespacing}{.5\linespacing}%
  {\normalfont\bfseries}}
\makeatother

\newcommand{\R}{\mathbb R}
\newcommand{\N}{\mathbb N}

\newcommand{\Z}{\mathbb Z}

\newcommand{\BEP}[1]{\textnormal{BEP}(#1)}

\newcommand{\BUCR}{\text{BUC}(\R)}

\newcommand{\arccosh}{\operatorname{arccosh}}

\newcommand{\gb}{g^{0}}
\newcommand{\Gd}{G^{\mathrm{dyn}}}
\newcommand{\gd}{g^{\mathrm{dyn}}}

\newcommand{\dd}{\hspace{2pt}\mathrm{d}}

\newcommand{\E}{\mathbb{E}}

\newcommand{\toP}{\overset{p}{\to}}

\newcommand{\p}{\mathbb{P}}

\newcommand{\U}{\mathcal U}
\newcommand{\xN}{x^{N}}
\newcommand{\yN}{y^{N}}
\newcommand{\zN}{z^{N}}

\begin{document}
\title[Hydrodynamic limits of ABEP and dynamic ABEP]{Hydrodynamic limits of the asymmetric Brownian energy process and its dynamic generalization}
\author{Frank Redig}
\author{Carel Wagenaar}
\date{\today}
\begin{abstract}
    We prove that the hydrodynamic limit of the asymmetric Brownian energy process (ABEP) is the viscous Burgers equation. We also prove the hydrodynamic limit of the dynamic ABEP which leads to a nonlinear equation for the macroscopic tail-energy potential. Both processes are related to the Brownian energy process (BEP) through invertible, non-local transformations involving discrete tail energies, which can be thought of as the microscopic equivalent of a transformation of Cole--Hopf type. Our main tool is a general micro-to-macro result showing that, under suitable conditions, such microscopic tail transformations induce corresponding transformations of the limiting density profiles. Combining this result with the fact that the hydrodynamic limit of BEP is the heat equation, we obtain the limiting profiles of ABEP and dynamic ABEP as explicit nonlinear transformations of solutions of the heat equation.
\end{abstract}
\maketitle

\section{Introduction} \label{sec:intro}
In the area of hydrodynamic limits one studies the emergence of partial differential equations from microscopic stochastic interacting particle systems. Under a suitable rescaling of space and time, the empirical density or energy field of such a stochastic microscopic system is shown to converge to a deterministic trajectory, which solves a partial differential equation (PDE). This PDE is then called the hydrodynamic limit of the system. Standard references on the hydrodynamic limits of interacting particle systems include \cite{Landim}, see also \cite{Demasi} and \cite{GR} for a duality-based approach to the hydrodynamic limit.
For symmetric systems such as the symmetric exclusion process (SEP), as well as for a class of zero-range processes, the rescaling of space and time is diffusive, and the hydrodynamic limit is a (possibly nonlinear) diffusion equation. For weakly asymmetric systems such as the weakly asymmetric exclusion process (WASEP), 
the hydrodynamic limit contains an additional (in general nonlinear) transport term. In \cite{Gar} the viscous Burgers equation is obtained as a hydrodynamic limit of the WASEP, using a microscopic equivalent of the Cole--Hopf transform which turns this nonlinear equation into the heat equation. In \cite{bertini} this approach is pushed further to prove that the density fluctuation field of the WASEP converges to the stochastic Burgers equation. See also \cite{CoGhMa} for the fluctuation limit of dynamic ASEP and \cite{CoShTs} for the fluctuation limit of $\text{ASEP}(q,j)$.
In this paper we add two new interacting particle systems for which the approach pioneered by G\"artner can be implemented, i.e., where on the microscopic level an analog  of a 
Cole--Hopf transformation exists. 

More precisely, in this paper, we study the hydrodynamic limit of the weakly asymmetric Brownian energy process and its dynamic generalization. The Brownian energy process, denoted by $\operatorname{BEP}(k)$, introduced in \cite{GKR}, is a symmetric diffusion process describing stochastic transport of a nonnegative continuous energy variable on a general graph with edge weights. In this paper, we restrict to  the case of the one-dimensional lattice $\Z$ with nearest-neighbor edges. The asymmetric version of the $\operatorname{BEP}(k)$, denoted by $\operatorname{ABEP}(\sigma,k)$, introduces a preferred direction of transport via an asymmetry parameter $\sigma$. Both processes arise as many-particle limits of inclusion-type particle systems, namely the symmetric inclusion process SIP$(k)$ and its asymmetric analog ASIP$(q,k)$ (see \cite{CGGRbep,CGRSsu11}).

A crucial feature of ABEP on the state space of configurations with finite total energy is its relation to BEP  via an invertible non-local transformation. This transformation is expressed in terms of the so-called discrete tail energies 
\[ 
    E_j(\vec y)=\sum_{i=j}^{\infty}y_i. 
\] 
More precisely, if $\vec Y(t)$ evolves according to ABEP, then the transformed process $\vec X(t)=g(\vec Y(t))$, with
\begin{equation}\label{ho1}
    g_j(\vec y) = \frac{e^{-2\sigma E_{j+1}(\vec y)} -e^{-2\sigma E_j(\vec y)}}{2\sigma}, 
\end{equation}
evolves according to BEP. Since this transformation is invertible on its natural range, the evolution of ABEP can be reconstructed from that of BEP. 

The dynamic asymmetric Brownian energy process, which we denote by $\operatorname{ABEP}_{\mathsf R}(\sigma,k,\rho)$, has similar properties. It depends on an additional right boundary parameter $\rho$ (explaining also the subindex $R$ in the acronym $\operatorname{ABEP}_{\mathsf R}(\sigma,k,\rho)$) and it arises as a scaling limit of dynamic ASIP introduced in \cite{Wa}. The latter is the inclusion version of dynamic ASEP, which was first introduced in \cite{BoCo} and its higher spin generalization in \cite{GrWaDyn}. The $\operatorname{ABEP}_{\mathsf R}(\sigma,k,\rho)$  is related to BEP through the invertible transformation 
\begin{equation}\label{ho2}
    g_j^{\mathrm{dyn}}(\vec z) = \frac{ \cosh\bigl(\sigma\rho+2\sigma E_j(\vec z)\bigr) - \cosh\bigl(\sigma\rho+2\sigma E_{j+1}(\vec z)\bigr) }{2\sigma}. 
\end{equation}
Both transformations have the form $f(E_j)-f(E_{j+1})$, and their inverses possess the same tail-difference structure. This observation is central to the micro-to-macro argument explained later on.

Although the processes ABEP and dynamic ABEP are considerably more complicated than their symmetric counterpart BEP, the nonlinear transformations \eqref{ho1}, resp.\ \eqref{ho2} make it possible to deduce several properties of ABEP and dynamic ABEP from the corresponding properties of BEP, such as reversibility and duality \cite{CGRSsu11, Wa}. 
This is in particular useful for the derivation of hydrodynamic limits, because the transformations \eqref{ho1}, \eqref{ho2} behave well in the scaling limit, i.e., converge to a nonlinear function of the empirical measure.

Using this approach, we show that the diffusive hydrodynamic behavior of ABEP and dynamic ABEP can be deduced from that of BEP (which is simple and given by the heat equation). In the diffusive setting, we rescale space by $N$ (i.e., the macro-point $x\in\R$ corresponds to the micro-point $\lfloor xN\rfloor)$ and rescale time with a factor $N^2$. With a configuration $\vec x\in\mathbb R_+^{\mathbb Z}$ we associate the empirical energy measure 
\[ 
    \Lambda_N(\vec x) = \frac1N\sum_{j\in\mathbb Z}x_j\delta_{j/N},
\] 
which is thought of as the macro-variable associated to the micro-configuration $\vec x\in\mathbb R_+^{\mathbb Z}$.
Then, under suitable assumptions on the initial configurations, including convergence of the empirical measure of these initial configurations to some absolutely continuous measure $\phi_0(x) \dd x$, the empirical measures of BEP converge to $\phi_t(x)\,\mathrm dx$, whose density solves the heat equation 
\[ 
    \partial_t\phi_t(x) = k\,\partial_{xx}\phi_t(x),
\] 
with initial datum $\phi_0$. 
\subsection{Micro and macro Cole--Hopf transformations}
Because we rescale time diffusively,   for the asymmetric systems, we will consider the weakly asymmetric regime
\[ 
    \sigma_N=\frac{\gamma}{2N}, 
\] 
where $\gamma\neq 0$ remains fixed as $N\to\infty$. We exploit the common structure of the microscopic transformations connecting BEP to ABEP and dynamic ABEP, as well as that of their inverses. These transformations take the form
\[ 
    h_j^N(\vec x) = \frac{N}{\gamma} \left[ \psi\left(\frac{\gamma}{N}E_j(\vec x)\right) - \psi\left(\frac{\gamma}{N}E_{j+1}(\vec x)\right) \right], 
\] 
for some smooth function $\psi$. Since $E_j(x) = E_{j+1}(x) + x_j$, the quantity 
$h_j^N(\vec x)$ is a discrete difference quotient. Formally, if the empirical measure of $\vec x^N$ converges to $\phi(x)\,\mathrm dx$, then 
\[ 
    \frac1N E_{\lfloor Nx\rfloor}(\vec x^{N}) \longrightarrow \int_x^\infty\phi(y)\,\mathrm dy. 
\] 
Consequently, the empirical measure of $h^N(\vec x^{N})$ is expected to converge to the density 
\[ 
    H(\phi)(x) = -\frac1\gamma \partial_x \psi\left( \gamma\int_x^\infty\phi(y)\,\mathrm dy \right) = \phi(x)\psi'\left( \gamma\int_x^\infty\phi(y)\,\mathrm dy \right). 
\] 
One of the contributions of this paper is a general theorem that makes this micro-to-macro correspondence rigorous. The result applies uniformly in time and allows us to transfer a known hydrodynamic limit 
(namely the heat equation for BEP) through a microscopic transformation towards a nonlinear equation such as
the Burgers equation and its dynamic variant.  

For the case of ABEP, the transformation $g$ converges in the weak-asymmetry regime to the invertible macroscopic map 
\[ 
    G(\mu)(x) = \mu(x) \exp\left( -\gamma\int_x^\infty\mu(y)\,\mathrm dy \right). \] 
We now start from initial configurations $\yN_0$ that converge to the macroscopic profile $\mu_0$. Then $g(\yN_t)$ follows BEP dynamics started from $g(\yN_0)$. By the theorem, $g(\yN_0)$ converges to
\[ 
    \phi_0=G(\mu_0). 
\] 
Thus in the hydrodynamic limit, $g(\yN_t)$ converges to $\phi_t$ which solves the heat equation with initial datum $\phi_0$. Then we use the theorem the other way around, which shows that $g^{-1}(g(\yN_t))$ converges to $$\mu_t:=G^{-1}(\phi_t).$$
Therefore, the function 
\[
    G(\mu_t) = \phi_t
\]
solves the heat equation. Since under appropriate conditions, the tail energy of a solution of the heat equation also solves the heat equation, we obtain that
\[ 
    J_t(x) :=1- \gamma \int_x^\infty G(\mu_t)(y)\dd y= \exp\left( -\gamma\int_x^\infty\mu_t(y)\,\mathrm dy \right) 
\] 
also satisfies the heat equation. This is a Cole--Hopf-type linearization and implies that $\mu_t$ solves the viscous Burgers equation 
\[ 
    \partial_t\mu_t = k\,\partial_{xx}\mu_t + k\gamma\,\partial_x\bigl(\mu_t^2\bigr). 
\] 
Alternatively, one could look at the tail potential $M_t=\int_x^\infty\mu_t(y)\dd y$, which satisfies the deterministic KPZ equation,
\begin{equation}\label{eq:detkpz}
    \partial_t M_t = k\partial_{xx} M_t - k\gamma (\partial_x M_t)^2.
\end{equation}

For the case of dynamic ABEP, the corresponding invertible macroscopic transformation reads
\[ 
    G^{\mathrm{dyn}}(\nu)(x) = \nu(x) \sinh\left( \rho+\gamma\int_x^\infty\nu(y)\,\mathrm dy \right). 
\] 
The natural macroscopic variable in this case is the tail-energy potential 
\[ 
    V_t(x) = \int_x^\infty\nu_t(y)\,\mathrm dy. 
\] 
After arguments similar to those used for ABEP, the quantity 
\[ 
    J_t^{\mathrm{dyn}}(x) = \cosh\bigl(\rho+\gamma V_t(x)\bigr) 
\] 
satisfies the heat equation. Substituting $J^\mathrm{dyn}_t$ into the heat equation gives that $V_t$ satisfies 
\begin{align}
    \partial_tV_t = k\,\partial_{xx}V_t + k\gamma \coth\bigl(\rho+\gamma V_t\bigr) \bigl(\partial_xV_t\bigr)^2. \label{eq:NCDKPZ}
\end{align}
This equation can be viewed as a nonlinear viscous Hamilton--Jacobi equation that is linearized by the hyperbolic transformation 
\[ 
    V\longmapsto\cosh(\rho+\gamma V). 
\] 
Differentiating the equation produces a non-local evolution equation for the density $\nu_t=-\partial_xV_t$, since the coefficient depends on the tail potential $V_t$. This explains why the tail-energy potential is the natural variable for the dynamic model. \\
When letting $\rho\to- \infty$ in this PDE, one obtains the deterministic KPZ equation \eqref{eq:detkpz}. In the limit $\rho\to \infty$, one obtains \eqref{eq:detkpz} with $\gamma$ replaced by $-\gamma$. Since this equation arises from dynamic ABEP, which is a process constructed in the realm of a representation of the non-compact Lie algebra $\mathfrak{su}(1,1)$, we refer to \eqref{eq:NCDKPZ} as the non-compact dynamic deterministic KPZ equation.

We thus derive that both ABEP and dynamic ABEP have similar linearization mechanisms. For ABEP, the exponential tail transformation reduces the viscous Burgers equation to the heat equation. For dynamic ABEP, the corresponding hyperbolic-cosine transformation reduces the tail-potential equation to the heat equation. Both mechanisms originate from deterministic transformations at the level of the microscopic configuration. 

The resulting correspondence can be summarized schematically by the commuting diagrams 
\[ 
    \begin{array}{ccc} \operatorname{ABEP} & \xleftrightarrow{\ g\ } & \operatorname{BEP}\\ \Big\downarrow\text{\scriptsize hydrodynamic limit} && \Big\downarrow\text{\scriptsize hydrodynamic limit}\\ \text{viscous Burgers equation} & \xleftrightarrow{\ G\ } & \text{heat equation} \end{array} 
\] 
and 
\[ 
    \begin{array}{ccc} \text{dynamic ABEP} & \xleftrightarrow{\ g^{\mathrm{dyn}}\ } & \operatorname{BEP}\\ \Big\downarrow\text{\scriptsize hydrodynamic limit} && \Big\downarrow\text{\scriptsize hydrodynamic limit}\\ \text{tail-potential equation} & \xleftrightarrow{\ G^{\mathrm{dyn}}\ } & \text{heat equation}. \end{array} 
\] 
\subsection{Main results}
Our main results can  be summarized as follows. 
\begin{enumerate} 
    \item We establish a general micro-to-macro theorem for non-local microscopic transformations expressed through discrete tail energies (Theorem \ref{thm:microtomacro}). 
    \item We use this theorem and the hydrodynamic limit of BEP to derive the hydrodynamic limit of weakly asymmetric ABEP under the scaling 
    \[
        \sigma_N = \frac{\gamma}{2N}.
    \]
    The limiting density is then found to be the unique solution of the viscous Burgers equation 
    \[ 
        \partial_t\mu = k\,\partial_{xx}\mu + k\gamma\,\partial_x(\mu^2), 
    \] 
    see Theorem \ref{thm:ABEP}.
    \item We derive the hydrodynamic limit of dynamic ABEP under the joint scaling 
    \[ 
        \sigma_N=\frac{\gamma}{2N}, \qquad \rho_N=\frac{2\rho N}{\gamma}. 
    \] 
    The limiting tail-energy potential solves 
    \[ 
        \partial_t V_t = k\,\partial_{xx}V_t + k\gamma\coth(\rho+\gamma V_t)(\partial_xV_t)^2,
    \] 
    see Theorem \ref{thm:dynABEP}. We will call this equation the ``non-compact dynamic deterministic KPZ equation''.
    \item We prove existence and uniqueness for the dynamic tail-potential equation in a natural class of bounded uniformly continuous functions (Theorem \ref{thm:uniquegenbur}). We also obtain an explicit solution formula and show that the evolution preserves monotonicity and the limits of the tail potential at spatial infinity, see Proposition \ref{prop:propertiesgenburger}.
\end{enumerate} 
The micro-to-macro theorem isolates the part of the argument that is specific to the non-local transformations from the part that concerns the hydrodynamic behavior of the underlying symmetric process. In particular, once the hydrodynamic limit of BEP is known, no direct replacement lemma or microscopic current computation is needed for ABEP or dynamic ABEP. The asymmetry is instead encoded in the transformation of the limiting profile. This provides a direct explanation for the appearance of the nonlinear macroscopic equations and may be useful for other stochastic systems related through tail-dependent transformations.

\subsection{Organization of the paper}
The remainder of the paper is organized as follows. In Section~\ref{sec:protra}, we introduce BEP, ABEP, and dynamic ABEP, together with the microscopic transformations connecting them. In Section~\ref{sec:mainresults}, we state the hydrodynamic-limit theorems. In Section~\ref{sec:microtomacro}, we prove the general micro-to-macro transformation theorem. In Section~\ref{sec:BEP}, we establish, or recall under the assumptions needed here, the hydrodynamic limit of BEP. Sections~\ref{sec:ABEP} and~\ref{sec:dynABEP} apply the transformation theorem to ABEP and dynamic ABEP, respectively. Finally, in Section~\ref{sec:uniexi}, we study existence, uniqueness, and qualitative properties of the nonlinear equation for the dynamic tail-energy potential.

\subsection{Some basic definitions and notation}
In this paper, we define BEP, ABEP and dynamic ABEP as diffusion processes on $\Z$ taking values in $$\R_+:=\{x\in \R: x\geq 0\}.$$ We will consider the state space of nonnegative sequences in $\ell^1(\Z)$, $$X:=\left\{ x \in \R_+^\Z: \sum_{j\in \Z} x_j < \infty\right\}.$$ Throughout this paper, BEP configurations will be denoted by $x$, ABEP configurations by $y$ and dynamic ABEP configurations by $z$.

Moreover, we define $$\big[\nabla_{i,j}f\big](x)=\big(\partial_{x_i} f-\partial_{x_{j}}f\big)(x).$$ 

Since we will look at the hydrodynamic limit, we will accelerate the process by $N^2$, where $N$ is the scaling factor. To simplify notation, when we have a Markov process $\{x_t: t\geq 0\}$ with generator $L$ and a sequence $\{\xN_0\}_{N\geq 1}$ of (random) initial states, we denote by $\xN_t$ the process started from $\xN_0$ at time $tN^2$. This rescaled  process 
$\{\xN_t: t\geq 0\}$ then has generator $N^2L$. 

Let $C_b(\R)$ denote the set of bounded continuous functions on $\R$.
Let $\mathcal{M}_+$ be the space of finite nonnegative Borel measures on $\R$ equipped with the weak topology, i.e.,  $\mu_N\to \mu$ if for all 
$f\in C_b(\R)$ we have 
\begin{equation}\label{weakl}
    \lim_{N\to\infty} \int_\R f \dd\mu_N = \int_\R f \dd\mu.
\end{equation}
Equivalently, \eqref{weakl} holds for all infinitely differentiable and compactly supported functions $f\in C_c^\infty(\R)$ together with the function $f=1$.
For a configuration $x\in X$, denote by $\Lambda_N(x)\in\mathcal{M}_+$ the empirical measure,
\[
    \Lambda_N(x)=\frac{1}{N} \sum_{j\in \Z} x_j \delta_{j/N},
\]
where $\delta_x$ is the Dirac measure. Then for $f\in C_b(\R)$ we have
\[
    \langle \Lambda_N(x),f\rangle:=\int_\R f\dd\Lambda_N(x) = \frac{1}{N} \sum_{j\in \Z} f(j/N)x_j.
\]  

By $\mu_N\toP \mu$ we denote weak convergence in probability as $N\to\infty$, i.e., for all $\varepsilon>0$
and all $f\in C_b(\R)  $
\[
\lim_{N\to\infty}\p\left(\left|\langle \mu_N-\mu,f\rangle  \right|>\varepsilon \right)=0.
\]
Let $T\geq 0$ be a finite time and for each $t\in[0,T]$, let $\mu_N(t)$, $N\geq 1$, be random measures in $\mathcal M_+$ and let
$\mu(t)$ be deterministic measures in $\mathcal M_+$. We say that
$\mu_N(t)\toP\mu(t)$ uniformly in $t\in[0,T]$ if for all $\varepsilon>0$ and all
$f\in C_b(\R)$
\begin{equation}\label{weakt}
    \lim_{N\to\infty}\p\left(\sup_{t\in[0,T]}\left|\langle\mu_N(t)-\mu(t),f\rangle\right|
    >\varepsilon\right)=0.
\end{equation}
If the family $\{\mu(t)\}_{t\in[0,T]}$ is uniformly bounded and tight, i.e.\
\[
    \sup_{t\in[0,T]}\mu(t)(\R)<\infty
    \qquad\text{and}\qquad
    \lim_{R\to\infty}\sup_{t\in[0,T]}\mu(t)\big(\{x:|x|>R\}\big)=0,
\]
then it suffices to verify \eqref{weakt} for $f=1$ and all $f\in C_c^\infty(\R)$. Both conditions hold in particular if
$t\mapsto\langle\mu(t),f\rangle$ is continuous on $[0,T]$ for $f=1$ and all
$f\in C_c^\infty(\R)$. In that case $t\mapsto\mu(t)$ is
continuous in $\mathcal M_+$, and \eqref{weakt} is equivalent to convergence in the Skorokhod path space
$D([0,T],\mathcal M_+)$.

For $\mu \in L^1(\R)$, the tail potential $M(x)=\int_x^\infty \mu(y) \dd y$ is bounded and uniformly continuous. Therefore, we introduce
\[
    \BUCR:=\{f\in C_b(\R)\colon f \text{ is uniformly continuous}\}.
\]
We define the heat kernel
\begin{align}
    K_t(x)= \frac{1}{\sqrt{4\pi k t}}e^{-\frac{x^2}{4kt}}. \label{eq:heatkern}
\end{align}
Then for $\phi_0 \in L^1(\R)$, the convolution
$$\phi_t(x) = (K_t * \phi_0)(x)  =\int_\R K_t(x-y)\phi_0(y) \dd y$$
is the unique solution in $C([0,T],L^1(\R))\cap C^{1,2}((0,T]\times \R)$ of the heat equation
\[
    \partial_t \phi_t = k \partial_{xx} \phi_t,\qquad (t,x)\in (0,T] \times \R,
\]
with initial condition $\phi_0$, thus
\[
    \lim_{t\downarrow 0} || \phi_t - \phi_0 ||_1 =0.
\]
Similarly if $\Phi_0\in \BUCR$, then $K_t*\Phi_0$ is the unique solution in $C([0,T],\BUCR)\cap C^{1,2}((0,T]\times \R)$ of the heat equation with initial condition $\Phi_0$, hence
\[
    \lim_{t\downarrow 0} || \Phi_t - \Phi_0 ||_\infty =0.
\]
Let us remark here already that if $\phi_0\in L^1(\R)$ and $\Phi_0(x)=\int_x^\infty \phi_0(y)\dd y$, then $\Phi_0 \in \BUCR$. In that case, if $\phi_t$ is the solution of the heat equation in $L^1(\R)$ as mentioned above, then $\Phi_t(x)=\int_x^\infty \phi_t(y)\dd y$ is the unique solution of the heat equation in $\BUCR$ as mentioned above. This will be proven in Lemma~\ref{lem:derheateqn} below.
\section{Processes and microscopic transformations}\label{sec:protra}
In this section, we define the relevant processes for this paper. We define their formal generators on the core of smooth cylinder functions on $\R^\Z_+$,
$$\mathcal{D}=\{f:\R^\Z_+ \to \R \colon f(x)=F(x_{j_1},\ldots,x_{j_m}),\ j_1,\ldots,j_m \in \Z,\  F\in C^\infty_c(\R^m_+) \}.$$

\subsection{Brownian energy process (BEP)}
Let $k>0$. The Brownian energy process BEP$(k)$ is the diffusion process on $X$ with generator given by
\begin{align}\label{bepgen}
	[L^{\mathrm{BEP}}_{k}f](x) = \sum_{j\in\Z}x_jx_{j+1} \big[(\nabla_{j,j+1})^2f\big](x) + k(x_{j+1}-x_{j}) \big[\nabla_{j,j+1} f\big](x).
\end{align}
The fact that this process is well-defined in the infinite-volume setting is proved e.g.  in \cite{zeg}. Alternatively, it can be shown by finite-volume approximations using duality with SIP, see \cite[Chapter 6]{GR}.
\subsection{Asymmetric Brownian energy process (ABEP)}
Let $\sigma\neq0$ and $k>0$, then ABEP$(\sigma,k)$ is the diffusion process on $X$ with generator given by 
\begin{equation}\label{abepgen}
\hspace{1cm}\big[L^{\mathrm{ABEP}}_{\sigma,k}f\big](y)= \sum_{j\in \Z}A_{j}^+ A_{j+1}^-  \big[(\nabla_{j,j+1})^2f\big](y) +\left( k[A_{j+1}^--A_{j}^+] - 2\sigma A_{j}^+A_{j+1}^-\right)\big[\nabla_{j,j+1} f\big](y),
\end{equation}
where
    \begin{align*}
        A_{j}^+&=\frac{1-e^{-2\sigma y_j}}{2\sigma} ,\quad  A_{j}^-= \frac{e^{2\sigma y_{j}}-1}{2\sigma}.
    \end{align*}
This process arises as many-particle and weak asymmetric limit of ASIP$(q,k)$, see \cite{CGRSsu11}. Since ABEP$(-\sigma,k)$ is the same as ABEP$(\sigma,k)$ with the order of sites reversed (i.e. $y_j\mapsto y_{-j}$), we can take $\sigma>0$ without loss of generality.

Let the partial energy $E_j(y)$ be defined as the total energy right of and including site $j$,
\[
	E_j(y)=\sum_{i=j}^\infty y_i,
\]
and the total energy by
\[
    E(y)=\sum_{j\in \Z} y_j.
\]
We will always consider configurations with finite total energy, i.e., $E(y)<\infty$. Since the processes we work with conserve the total energy, this condition is also satisfied at later times, and it implies that the tail energies $E_j(y)$ are all well-defined. 
Define $X_\sigma:=\{ x\in X: E(x) < \frac{1}{2\sigma}\}$. Then the mapping $g\colon X\to X_\sigma$ defined by
\begin{align}
	g_j(y)=\frac{e^{-2\sigma E_{j+1}(y)}- e^{-2\sigma E_{j}(y)} }{2\sigma },
\end{align}
is a bijective transformation between BEP and ABEP. More precisely, from the  generator equality 
\begin{align}
	\big[L^\mathsf{BEP}_kf\big](g(y))= \big[L^\mathsf{ABEP}_{\sigma,k}(f\circ g)\big](y),\label{eq:ABEPBEPtransformation}
\end{align}
one derives that if $Y(t)$ follows ABEP dynamics, $g(Y(t))$ follows BEP dynamics (see \cite{CGRSsu11}).
The codomain $X_\sigma$ arises via the computation
\[
    E(g(y))=\sum_{j=-\infty}^\infty g_j(y)= \frac{1-e^{-2\sigma E(y)} }{2\sigma } < \frac{1}{2\sigma }.
\]
The inverse $g^{-1}: X_\sigma\to X$ has components given by
\begin{align}
    g_j^{-1}(x)= \frac{\log\left(1- 2\sigma E_{j+1}(x)\right) -\log\left(1- 2\sigma E_{j}(x)\right) }{2\sigma} \label{eq:ginverse}.
\end{align}
This follows from
\[
    2\sigma E_j(g(y)) = 2\sigma \sum_{i=j}^\infty g_i(y)= 1-  e^{-2\sigma E_{j}(y)},
\]
combined with
\[
      E_j(y) = -\frac{\log(1-2\sigma E_j(g(y)))}{2\sigma},
\]
and $y_j = E_{j}(y)-E_{j+1}(y)$.

Notice that the transformation $g^{-1}$ can be used to \textit{define} the process with formal generator
\eqref{abepgen} on the space $X$ of configurations with finite total energy, by taking  $Y(t)=g^{-1}(X(t))$, where $X(t)$ is a BEP process. In this way we circumvent the (non-trivial) problem of defining the ABEP process in infinite volume directly from its generator \eqref{abepgen}, which would possibly define this process on a larger state space also including configurations with infinite total energy. However, we do not need this here since we will always consider the ABEP on configuration with finite total energy.
\subsection{Dynamic ABEP}
    Let $k,\rho >0$ and $\sigma\neq0$, where $\rho$ can be interpreted as a `right boundary value'. Then we define ABEP$_\mathsf{R}(\sigma,k,\rho)$, the right version of the dynamic ABEP, as the diffusion process on $X$ with generator given by
	\begin{equation}\label{dabepgen}
	[L^{\mathrm{ABEP}_\mathsf{R}}_{\sigma,k,\rho}f](z)= \sum_{j\in \Z} C_{j}^+ C_{j+1}^-  \big[(\nabla_{j,j+1})^2f\big](z) +\left( k[C_{j+1}^--C_{j}^+] + C_{j}^+C_{j+1}^-D_j\right)\big[\nabla_{j,j+1} f\big](z),
	\end{equation}
	where
    \begin{align*}
        C_{j}^+&=\frac{\sinh(\sigma z_{j})}{\sigma}\frac{\sinh\big(\sigma [\rho+2E_{j+1}(z)+z_{j}]\big)}{\sinh\big(\sigma [\rho+2E_{j+1}(z)]\big)} ,   \\
        C_{j+1}^-&= \frac{\sinh(\sigma z_{j+1})}{\sigma} \frac{\sinh\big(\sigma [\rho+2E_{j+1}(z)-z_{j+1}]\big)}{\sinh\big(\sigma[\rho+2E_{j+1}(z)]\big)}, \\
        D_j&= \frac{2\sigma\cosh\big(\sigma[\rho+2E_{j+1}(z)]\big)}{\sinh(\sigma[\rho+2E_{j+1}(z)])}.
    \end{align*}
\begin{remark}\* \label{rem:dynABEP}
	\begin{enumerate}[label=(\roman*)]
        \item Since $C_j^+,C^-_{j+1}$ and $D_j$ are all invariant under the transformation $\sigma\mapsto -\sigma$, ABEP$_\mathsf{R}(-\sigma,k,\rho)$ is the same as ABEP$_\mathsf{R}(\sigma,k,\rho)$. Therefore, we can take $\sigma >0$ without loss of generality.
        \item In general $\rho\in \R$ is allowed, as long as the coefficient $C_j^+C^-_{j+1}$ in front of the $\nabla_{j,j+1}^2$ remains nonnegative and does not blow up. Taking $\rho>0$ is a natural choice for this. If one would take $\rho <0$, this can be done by for instance limiting the state space to states $z\in X$ such that $2E(z) < |\rho|$, since then $\rho+2E_j(z)<0$ for all $j\in\Z$. In this paper, we will choose the first option and take $\rho >0$. \label{it:rho<0}
		\item In the limit $\rho\to-\infty$ one obtains ABEP$(\sigma,k)$ and when $\rho\to\infty$ one gets ABEP$(-\sigma,k)$.
		\item ABEP, resp.\ dynamic ABEP arises as a weak asymmetric limit of ASIP, resp.\ dynamic ASIP, see \cite{CGRSsu11}, \cite{Wa}.
        \item One can reverse the order of sites ($x_j \mapsto x_{-j}$). Then $\rho$ becomes a `left' boundary value and one gets the left version of dynamic ABEP.
	\end{enumerate}
\end{remark}

The explicit form of the generator is a bit complicated and often inconvenient to work with. Therefore, we will exploit that for dynamic ABEP, there is also a non-local transformation connecting it to BEP (\cite{Wa}):
\begin{align}
	\big[L^\mathrm{BEP}_kf\big](g^{\mathrm{dyn}}(z))= \big[L^{\mathrm{ABEP}_\mathsf{R}}_{\sigma,k,\rho}(f\circ g^{\mathrm{dyn}})\big](z),\label{eq:LABEPBEPtransformation}
\end{align}   
where $\gd:X \to X$ has components given by
\begin{align*}
	\gd_j(z)= \frac{\cosh(\sigma\rho+2\sigma E_{j}(z))-\cosh(\sigma\rho + 2\sigma E_{j+1}(z))}{2\sigma}.
\end{align*}
Consequently, if $Z(t)$ follows ABEP$_\mathsf{R}$ dynamics, $\gd(Z(t))$ follows BEP dynamics. This non-local transformation generalizes the one between ABEP and BEP. Indeed, multiplying $\gd$ by $-2e^{\sigma\rho}$ and letting $\rho\to-\infty$ gives the desired transformation. Note that $\gd(z)$ has finite energy,
\[
    E(\gd(z))=\sum_{j=-\infty}^\infty \gd_j(z) = \frac{\cosh(\sigma \rho + 2\sigma E(z))-\cosh(\sigma \rho )}{2\sigma}.
\]
Moreover, $\gd \colon X \to X$ is invertible, with inverse 
\[
    (\gd)^{-1}_j(x) = \frac{\arccosh(\cosh(\sigma \rho) + 2\sigma E_j(x)) - \arccosh(\cosh(\sigma \rho) + 2\sigma E_{j+1}(x)) }{2\sigma},
\]
where $\arccosh\colon [1,\infty)\to[0,\infty)$ is the positive branch of the inverse of $\cosh$. Here, $\sigma\rho >0$ ensures that $\arccosh(\cosh(\sigma\rho))=\sigma\rho$.

Similarly to the case of ABEP, we can define the dynamic ABEP process by taking  $Z(t)=(\gd)^{-1}(X(t))$, where $X(t)$ is a BEP process. 

\section{Main results}\label{sec:mainresults}
In this section we will state our main theorems: the hydrodynamic limits of ABEP and dynamic ABEP. We take the hydrodynamic limit in a weak asymmetry regime, i.e. $\sigma_N = \gamma/2N$. 
\begin{theorem}[ABEP] \label{thm:ABEP}
    Assume that the sequence of empirical measures corresponding to (random) initial configurations $\{\yN_0\}_{N\geq1}$ converges in probability to a deterministic measure which is absolutely continuous with respect to the Lebesgue measure,
    \begin{align}
        \Lambda_N(\yN_0) \toP \mu_0(y)\dd y \qquad \text{in }\mathcal{M}_+. \label{eq:ABEPinitialconvergence} 
    \end{align}
    Moreover, assume that the second moment of $(\yN_0)_j$ is uniformly bounded in $j$ and $N$, i.e. there exists a constant $C_\mathrm{A}$ such that 
        \begin{align}
            \sup_{j\in \Z, N\in \N} \E\bigg[\big(\yN_0\big)_j^2\bigg] \leq C_{\mathrm{A}}, \label{eq:ABEPboundquadexpectation}
        \end{align}
        and the total energy grows at most linearly in $N$, i.e. there exists a constant $E_\mathrm{A}$ such that almost surely for all $N$ we have
        \begin{align}
            \frac{1}{N}\sum_{j\in \Z} (\yN_0)_j \leq E_{\mathrm{A}}. \label{eq:ABEPboundtotalenergy}
        \end{align}
        Let $\gamma\neq0$ and define $\sigma_N=\gamma/2N$. Then for each time $T>0$, the empirical measure of the accelerated ABEP$(\sigma_N,k)$ process started at $\yN_0$, converges in probability to the measure $\mu_t(y)\dd y$, uniformly in $t\in [0,T]$,
    \[
        \Lambda_N(\yN_t) \toP \mu_t(y)\dd y \qquad \text{in }\mathcal{M}_+ \text{ uniformly in } t\in [0,T],
    \]
    where $\mu_t$ is the unique solution in $C([0,T],L^1(\R))\cap C^{1,2}((0,T]\times \R)$ of the viscous Burgers equation
    \begin{align*}
        	&\partial_{t}\mu_t = k \partial_{yy}\mu_t + k\gamma  \partial_{y}(\mu_t^2),\qquad (t,y) \in (0,T] \times \R,
    \end{align*}
    satisfying 
    \[
        \lim_{t\downarrow 0} || \mu_t - \mu_0||_1 =0.
    \]
\end{theorem}
     \begin{remark}\* \label{rem:muintegrable}
         \begin{enumerate}[label=(\roman*)]
            \item Proving the theorem for $\gamma>0$ is enough, since the case $\gamma<0$ follows by the spatial reflection $y_j\mapsto y_{-j}$. Indeed, this reflection changes $\operatorname{ABEP}(\sigma_N,k)$ to $\operatorname{ABEP}(-\sigma_N,k)$. The assumptions on the second moments and total energy are invariant under this reflection. Moreover, if the empirical measures of the initial configurations converge to $\mu_0(y)\,\mathrm{d}y$, then those of the reflected configurations converge to $\mu_0(-y)\,\mathrm{d}y$. At the macroscopic level, if $\mu_t$ solves \[ \partial_t\mu_t = k\partial_{yy}\mu_t + k\gamma\partial_y(\mu_t^2), \] then $\widetilde\mu_t(y)=\mu_t(-y)$ solves the same equation with $\gamma$ replaced by $-\gamma$. Thus, changing the sign of $\gamma$ only reverses the direction of the asymmetry.
            \item Setting $\tilde{\mu}_t= 2|\gamma|\mu_t$, $\tau = kt$, and using the chain rule gives a rescaled viscous Burgers equation of the form $$\partial_\tau \tilde{\mu}_\tau = \partial_{yy}\tilde{\mu}_\tau + \text{sign}(\gamma) \tilde{\mu}_\tau(\partial_y \tilde{\mu}_\tau).$$
             \item Note that $\Lambda_N(\yN_0)\toP \mu_0(y)\dd y$ in combination with assumption \eqref{eq:ABEPboundtotalenergy} imply that $\mu_0 \in L^1(\R)$ with $$\int_\R \mu_0(y) \dd y  \leq E_\mathrm{A},$$
         by using weak convergence with the test function $f=1$. This also explains the $L^1(\R)$ appearing in the solution class $C([0,T],L^1(\R))$.
         \item Assumption \eqref{eq:ABEPboundquadexpectation} implies $\mu_0 \in L^\infty(\R)$. Indeed, for any nonnegative $f\in C_b(\R)$ with compact support we can extract a subsequence $y^{N_m}_0$ such that 
         \[
                \lim_{m\to\infty}\langle \Lambda_{N_m}(y^{N_m}_0),f\rangle = \int_\R f(x)\mu_0(x)\dd x \quad \text{a.s.}
         \]
         Then, by combining Fatou's lemma with the Cauchy-Schwarz inequality, we obtain
         \begin{align*}
             \int_\R f(x)\mu_0(x) \dd x &\leq \liminf_{m\to\infty} \E\left[\langle \Lambda_{N_m}(y^{N_m}_0),f\rangle \right]\\ &=\liminf_{m\to\infty} \frac1{N_m}\sum_{j\in \Z}  f(j/N_m) \E[(y^{N_m}_0)_j]\\
             &\leq \sqrt{C_\mathrm A} \int_\R f(x)\dd x,
         \end{align*}
         where in the last step we used that for 
         $f\in\ C_b(\R)$ with compact support, the Riemann sums $\frac1N\sum_{j\in \Z}  f(j/N)$ converge to the integral. Therefore, $||\mu_0||_\infty \leq \sqrt{C_\mathrm A}$.
        \end{enumerate}
     \end{remark}
     In the hydrodynamic limit of dynamic ABEP we let $\sigma_N\to 0$. However, we do not want $\sigma_N\rho$ to vanish. Therefore, we take $\rho_N = 2\rho N/\gamma$ so that $\sigma_N\rho_N=\rho$.
\begin{theorem}[Dynamic ABEP]\label{thm:dynABEP}
    Assume that the sequence of empirical measures corresponding to (random) initial configurations $\{\zN_0\}_{N\geq1}$ converges in probability to a deterministic measure which is absolutely continuous with respect to the Lebesgue measure,
    \[
        \Lambda_N(\zN_0) \toP \nu_0(z)\dd z \qquad \text{in }\mathcal{M}_+.
    \]
    Moreover, assume that the second moment of $(\zN_0)_j$ is uniformly bounded in $j$ and $N$, i.e. there exists a constant $C_{\mathrm{D}}$ such that 
        \begin{align}
            \sup_{j\in \Z, N\in \N} \E\bigg[\Big(\zN_0\Big)_j^2\bigg] \leq C_\mathrm{D}, \label{eq:dynABEPboundquadexpectation}
        \end{align}
        and the total energy grows at most linearly in $N$, i.e. there exists a constant $E_\mathrm{D}$ such that almost surely for all $N$ we have
        \begin{align}
            \frac{1}{N}\sum_{j\in \Z} (\zN_0)_j \leq E_{\mathrm{D}}. \label{eq:dynABEPboundtotalenergy}
        \end{align}
    Let $\gamma,\rho >0$ and define $\sigma_N=\gamma/2N$ and $\rho_N=\rho/\sigma_N$. Then for each time $T>0$, the empirical measure of the accelerated ABEP$_\mathsf{R}(\sigma_N,k,\rho_N)$ process started at $\zN_0$, converges in probability to the measure $\nu_t(y)\dd y$, uniformly in $t\in[0,T]$,
    \[
        \Lambda_N(\zN_{t}) \toP \nu_t(z)\dd z \qquad \text{in }\mathcal{M}_+ \text{ uniformly in } t\in [0,T].
    \]
    Here $V_t(z)=\int_z^\infty \nu_t(\tilde{z})\dd \tilde{z}$ is the unique solution in $C([0,T],\BUCR)\cap C^{1,2}((0,T]\times \R)$ of
        the non-compact dynamic deterministic KPZ equation
    \begin{align*}
        	&\partial_{t} V_t(z) = k \partial_{zz}V_t(z) +   k\gamma\coth(\rho + \gamma V_t(z)) \big(\partial_{z}V_t(z)\big)^2
    \end{align*}
    satisfying 
    \[
            \lim_{t\downarrow 0} || V_t - V_0||_\infty =0.
    \]
\end{theorem}
\begin{remark}\*
\begin{enumerate}[label=(\roman*)]
    \item One could take $\gamma <0$. However, then $\rho_N = 2\rho N/\gamma<0$ as well, and one should make sure that the rate coefficients $C_j^+C^-_{j+1}$ of dynamic ABEP are nonnegative and do not blow up. As explained in Remark~\ref{rem:dynABEP}\ref{it:rho<0}, this can be done by requiring that $2E(z) < |\rho_N|$ which is equivalent to $$\rho+\frac{\gamma}{N}E(z) > 0.$$ Note that this implies that  $\rho + \gamma V_t(z) >0$ as well, hence $\coth(\rho + \gamma V_t(z))$ is well-defined in the resulting PDE.
    \item Setting  $\tilde{V}_{t}(z)=\gamma V_t(z)+\rho$ and $\tau=kt$ gives the rescaled PDE 
    \[
    \partial_\tau \tilde{V}_\tau(z) = \partial_{zz} \tilde{V}_\tau(z) + \coth(\tilde{V}_\tau(z))\big(\partial_z \tilde{V}_\tau(z)\big)^2.
    \]
    \item Similar to Remark~\ref{rem:muintegrable}, we have $\nu_0 \in L^1(\R) \cap L^\infty(\R)$ with $||\nu_0||_1 \leq E_\mathrm D$ and $||\nu_0||_\infty \leq \sqrt{C_\mathrm D}$. This in turn implies $V_0 \in \BUCR$, which leads to the solution class $C([0,T],\BUCR)$ of the PDE. 
\end{enumerate}
\end{remark}
    \section{A general micro-to-macro transformation theorem} \label{sec:microtomacro}
    
    In this section, we prove that the microscopic transformations $g,g^{-1},\gd,(\gd)^{-1}$ all have a macroscopic limit. The key observation here is that all these transformations are of the form $h=(h_j^N)_{j\in \Z}$, where the $j$-th component is given by
    \[
        h_j^N(x) = \frac{N}{\gamma}\left(\psi\left(\frac{\gamma}{N}E_{j}(x)\right)-\psi\left(\frac{\gamma}{N}E_{j+1}(x)\right)\right),
    \]
    where $\psi$ is a smooth function on some compact interval. Let us first do a formal computation for the microscopic transformation $h^N$. Assume that $\Lambda_N(\yN)\toP \mu(y)\dd y$ and the microscopic site $j_N$ converges to the macroscopic site $y$, i.e. $j_N/N \to y$. Let $f$ be a test function, then
        \begin{align}
        	\langle \Lambda_N(h(\yN)) , f \rangle  = \frac{1}{N}\sum_{j=-\infty}^\infty f(j/N)\frac{N}{\gamma}\left(\psi\left(\frac{\gamma}{N} E_{j}(\yN)\right)- \psi\left( \frac{\gamma}{N}E_{j+1}(\yN)\right)\right) \label{eq:ABEPtest}
        \end{align}
        Since $j_N/N\to y$, we have
        \begin{align}
            \frac{1}{N}E_{j_N}(\yN) \to \int_{y}^\infty \mu(x)\dd x, \label{eq:formalconvpartialenergy}
        \end{align}
        thus
        \begin{align*}
        \psi\left(\frac{\gamma}{N} E_{j_N}(y^N)\right) \to \psi\left(\gamma \int_{y}^\infty \mu(x)\dd x\right).
        \end{align*}
        Since $E_{j+1}(y)=E_{j}(y) - y_j$, the coefficient $h_j^N$ is a discrete difference quotient. Therefore,
        \[
                N \left[\psi\left(\frac{\gamma}{N} E_{j}(\yN)\right)- \psi\left( \frac{\gamma}{N}E_{j}(\yN)-\frac{\gamma}{N}\yN_j\right)\right] \to - \partial_y \left[ \psi\left( \gamma \int_y^\infty \mu(x)\dd x\right) \right].
        \]
        Hence \eqref{eq:ABEPtest} converges to
        \begin{align*}
        	\int_\R f(y)H(\mu)(y) \dd y,
        \end{align*}
        where $H$ is the corresponding macroscopic transformation of $h^N$ given by
        \begin{align*}
            H(\mu)(y) &= -\frac{1}{\gamma} \partial_y\left[\psi\left(\gamma \int_{y}^\infty \mu(x)\dd x\right)\right] \\
            &=\mu(y)\psi'\left(\gamma \int_{y}^\infty \mu(x)\dd x\right).
        \end{align*}
        Note that if $\mu\in L^1(\R)$, the first representation of $H$ including the derivative is well-defined by the Lebesgue differentiation theorem. In Table~\ref{tab:micro-macro-transformations}, an overview of the macroscopic transformations of the four different microscopic transformations from BEP to (dynamic) ABEP and back can be found.
        
        \begin{table}[ht] \centering \renewcommand{\arraystretch}{1.8} \begin{tabular}{c|c|c} \hline $h_j^N$ & $\psi(a)$ & $H(u)(x)$ \\ \hline $g_j$ & $-e^{-a}$ & $\displaystyle G(u)(x)=u(x)e^{-\gamma\int_x^\infty u(y)\,\mathrm{d}y}$ \\ $g_j^{-1}$ & $\displaystyle -\log(1-a)$ & $\displaystyle G^{-1}(u)(x) =\frac{u(x)} {1-\gamma\int_x^\infty u(y)\,\mathrm{d}y}$ \\ $g_j^{\mathrm{dyn}}$ & $\cosh(\rho+a)$ & $\displaystyle \Gd(u)(x)= u(x)\sinh\left( \rho+\gamma\int_x^\infty u(y)\,\mathrm{d}y \right)$ \\ $\left(g^{\mathrm{dyn}}\right)^{-1}_j$ & $\displaystyle \arccosh\!\left(\cosh(\rho)+a\right)$ & $\displaystyle (\Gd)^{-1}(u)(x)= \frac{u(x)} {\sqrt{ \left( \cosh(\rho) +\gamma\int_x^\infty u(y)\,\mathrm{d}y \right)^2-1 }}$ \\ \hline \end{tabular} \caption{Microscopic transformations and their macroscopic counterparts.} \label{tab:micro-macro-transformations} \end{table}  
        To make this computation rigorous, we first prove \eqref{eq:formalconvpartialenergy}.
        \begin{lemma}\label{lem:convpartialenergy}
            Let $T\geq 0$ and assume that for each $t\in [0,T]$, $\{\mu^N_t\}_{N\geq 1}$ is a sequence of random measures which converges in probability to the random measure $\mu_t$ uniformly in $t\in [0,T]$, where $\mu_t$ satisfies for each $x\in \R$  \begin{align}\lim\limits_{\delta \downarrow 0}\sup\limits_{t\in[0,T]} \mu_t([x-\delta,x+\delta]) =0 \label{eq:uniformnoatoms}\end{align} in probability. Moreover, let $M\in \N$ and assume that $j_N\in [-MN,MN]$ is a sequence of integers such that $j_N/N\to y\in [-M,M]$. Then for all $\varepsilon >0$ we have
            \begin{align*}
                \lim_{N\to\infty} \p\Big( \sup_{t\in[0,T]} \big|\mu^N_t([j_N/N,\infty))- \mu_t([y,\infty))\big| > \varepsilon\Big) =0.
            \end{align*}
            Additionally, if $a_N(j)= \p\Big( \sup_{t\in [0,T]}\big|\mu_t^N([j/N,\infty))- \mu_t([j/N,\infty))\big| > \varepsilon\Big)$, then 
            \begin{align*}
                \lim_{N\to\infty}\frac{1}{N}\sum_{j=-NM}^{NM}a_{N}(j) =0.
            \end{align*}
        \end{lemma}
        \begin{proof}
        Let us first show that for each $x\in\R$, we have 
        \begin{align}
            \mu_t^N\big([x,\infty) \big)\toP\mu_t\big([x,\infty)\big), \label{eq:muNt}
        \end{align} 
        uniformly in $t\in [0,T]$. Indeed, let $f_n^-, f_n^+\in C_b(\R)$ such that $0\leq f_n^-,f_n^+ \leq 1$ and 
        \begin{align*}
            f_n^{+}(y) &= \begin{cases} 0\qquad \text{if } y\leq x-\tfrac1n,\\ 1 \qquad\text{if } y \geq x, \end{cases} \\
             f_n^{-}(y) &= \begin{cases} 0\qquad \text{if } y\leq x,\\ 1 \qquad\text{if } y \geq x +\tfrac1n. \end{cases}   
        \end{align*}
        Then $f_n^- \leq 1_{[x,\infty)} \leq f_n^+$. For the upper bound,
        \begin{align*}
            \mu^N_t([x,\infty)) \leq \int f_n^+ \dd \mu^N_t = \int f_n^+ \dd \mu_t + \left(\int f_n^+ \dd \mu^N_t - \int f_n^+ \dd \mu_t \right) \\
            \leq \mu_t([x-\tfrac1n,x)) + \mu_t([x,\infty)) + \left(\int f_n^+ \dd \mu^N_t - \int f_n^+ \dd \mu_t \right).
        \end{align*}
        Therefore,
        \[
            \mu^N_t([x,\infty)) - \mu_t([x,\infty)) \leq \mu_t([x-\tfrac1n,x)) + \left(\int f_n^+ \dd \mu^N_t - \int f_n^+ \dd \mu_t \right).
        \]
        The right hand side goes to $0$ uniformly in $t\in[0,T]$ by \eqref{eq:uniformnoatoms} and $\mu_t^N \toP \mu_t$ uniformly in $t\in[0,T]$. Now proceed in a similar way with $f_n^-$ for the lower bound, proving \eqref{eq:muNt}. 

        Now let $\varepsilon,\varepsilon'>0$. By \eqref{eq:uniformnoatoms}, we can find a $\delta>0$ such that 
        \begin{equation}\label{pong}
            \p\left(\sup\limits_{t\in[0,T]} \mu_t([y-\delta,y+\delta]) > \varepsilon/2\right) < \varepsilon'/2.
        \end{equation}
        Then, take $N'$ such that for all $N\geq N'$ we have $j_N/N \in [y-\delta,y+\delta]$ and
        \begin{equation}\label{ping}
            \p\left(\sup_{t\in[0,T]} |\mu_t^N([y\pm \delta,\infty))-\mu_t([y\pm\delta,\infty))|>\varepsilon/2\right) < \varepsilon'/2.
        \end{equation}
        By monotonicity of $\mu^N_t$, we have
        \begin{align*}
            \mu_t^N\big([y+\delta,\infty)\big) \leq \mu^N_t \big([j_N/N,\infty)\big) \leq \mu_t^N\big([y-\delta,\infty)\big) 
        \end{align*}
        Therefore,
        \begin{align*}
             \mu_t^N([j_N/N,\infty)) - \mu_t([y,\infty))&\leq \mu_t^N([y-\delta,\infty)) - \mu_t([y-\delta,\infty)) +\mu_t([y-\delta,y+\delta]),
        \end{align*}
        where we used that 
        \[
        \mu_t([y,\infty))= \mu_t([y-\delta,\infty)) -\mu_t([y-\delta,y)) \geq \mu_t([y-\delta,\infty))-\mu_t([y-\delta,y+\delta]).
        \]
        Similarly,
        \begin{equation*}
        \mu_t^N([j_N/N,\infty)) - \mu_t([y,\infty)) 
        \geq \mu^N_t ([y+\delta, \infty)) - \mu_t ([y+\delta, \infty)) - \mu_t([y-\delta,y+\delta]).
        \end{equation*}
        Combining this with \eqref{pong}, \eqref{ping} gives $\mu_t^N \big([j_N/N,\infty)\big)\toP \mu_t\big([y,\infty)\big)$ uniformly in $t\in[0,T]$.
        
        For the second part, we argue by contradiction. Suppose that $$\frac{1}{N}\sum_{j=-NM}^{NM} a_{N}(j)$$ does not converge to $0$. Then there exists a $\delta>0$ and a subsequence $N_k$ such that $$ \Big|\frac{1}{N_k}\sum_{j=-N_kM}^{N_kM} a_{N_k}(j)\Big| \geq \delta$$ for all $k\in \N$. Therefore, we can find for each $k$ a $j_k\in [-N_kM,N_kM]$ such that $a_{N_k}(j_k) \geq \delta/(2M+1)$. Define $y_k=j_k/N_k$, then $y_k\in[-M,M]$, thus by compactness we can pass to a further subsequence so that $y_k \to y\in [-M,M]$. However, then $a_{N_k}(j_{k})\to 0$ by the first part of the lemma and \eqref{eq:uniformnoatoms}, giving a contradiction.
        \end{proof}
        \begin{remark}
            In our setting, $\mu_t= \phi_t(x) \dd x$ is absolutely continuous with respect to the Lebesgue measure. Then $$\lim\limits_{\delta\downarrow 0}\sup\limits_{t\in[0,T]}\mu_t([x-\delta,x+\delta]) =0$$ is automatically satisfied when $T=0$ and also when $T>0$ whenever $\sup\limits_{t\in [0,T]} ||\phi_t||_\infty < \infty$.
        \end{remark}
        Now we are ready to prove the micro to macro result.
        \begin{theorem} \label{thm:microtomacro}
            Let $T\geq 0$. For every $N$, let $t\mapsto \xN_t$ be an $\mathbb R_+^{\mathbb Z}$-valued random process such that 
            \[
                \sup_{t\in[0,T]} \left|\langle \Lambda_N(\xN_t),f \rangle -\int_\R f(x)\phi_t(x)\dd x\right| \toP 0 ,
            \]
            for all $f\in C_c^\infty(\R)$ and $f=1$, where the deterministic measurable functions $\{\phi_t\}_{t\in[0,T]}$ satisfy
            \begin{align} 
            \sup_{t\in [0,T]} ||\phi_t||_\infty < \infty,\qquad
            \lim_{R\to\infty} \sup_{t\in [0,T]} \int_{|x|>R} \phi_t(x)\dd x =0. \label{eq:ubddtight}
            \end{align}
            Moreover, assume that there exists a constant $E$ such that for all $t\in[0,T]$ and all $N$ we almost surely have
            \begin{align}
                \frac{1}{N} \sum_{j\in \Z} (\xN_t)_j &\leq E.\label{eq:Eass}
            \end{align}
            Fix $\gamma >0$, a non-decreasing $\psi\in C^2([0,\gamma E])$ and define $h^N=(h_j^N)_{j\in\Z}$ and $H$ by
            \begin{align*}
                h_j^N(x)&= \frac{N}{\gamma}\left(\psi\left(\frac{\gamma}{N}E_j(x)\right)-\psi\left(\frac{\gamma}{N}E_{j+1}(x)\right)\right),\\
                H(\phi_t)(x) &= -\frac{1}{\gamma} \partial_x\left[\psi\left(\gamma \int_{x}^\infty \phi_t(y)\dd y\right)\right] \\
                        &=\phi_t(x)\psi'\left(\gamma \int_{x}^\infty \phi_t(y)\dd y\right).
            \end{align*}
            Then 
            \[
                \Lambda_N(h^N(\xN_t)) \toP H(\phi_t)(x)\dd x \quad \text{in }\mathcal{M}_+ \text{ uniformly in }t\in [0,T].
            \]
        \end{theorem}  
        \begin{proof}
            Note that $\phi_t(x)\geq 0$ and, similar to Remark~\ref{rem:muintegrable}, we have $\int \phi_t(x)\dd x \leq E$. Pick $\varepsilon>0$. Because the measures $\{ \phi_t(x)\dd x\}_{t\in[0,T]}$ are uniformly bounded and tight and
            \[
                |H(\phi_t)| \leq |\phi_t| ||\psi'||_\infty,
            \]
            the measures $\{ H(\phi_t)(x)\dd x\}_{t\in[0,T]}$ inherit these properties. Therefore,
            it is enough to show that 
            \[
                \lim_{N\to\infty} \p\bigg[\sup_{t\in[0,T]}\Big|\frac{1}{N}\sum_{j\in \Z}h_j(\xN_t)f(j/N) - \int_\R H(\phi_t)(x)f(x)\dd x \Big|>\varepsilon \bigg] =0
            \]
            for $f=1$ and $f\in C_c^\infty(\R)$. For $f=1$, note that by telescoping
            \[
                    \frac{1}{N}\sum_{j\in \Z} h_j(\xN_t) =  \frac{1}{\gamma} \psi \left( \frac{\gamma}{N}  E(\xN_t)\right) - \frac{1}{\gamma}\psi(0),
            \]
            where $\frac{1}{N}E(\xN_t) \leq E$ by assumption \eqref{eq:Eass}. For the other term, we have
            \[
                    \int_\R H(\phi_t)(x) \dd x = \frac{1}{\gamma} \psi\left(\gamma \int_\R \phi_t(y)\dd y\right) - \frac{1}{\gamma}\psi(0).
            \]
            Since $|\psi(x)-\psi(y)|\leq |x-y|\cdot ||\psi'||_\infty$,
            \[
                \left|\psi\left(\frac{\gamma}{ N}  E(\xN_t)\right) - \psi\left(\gamma \int_\R \phi_t(y)\dd y\right)\right| \leq \gamma ||\psi'||_\infty  \left|\frac{1}{ N}  E(\xN_t) - \int_\R \phi_t(y)\dd y\right|.
            \]
            Therefore,
            \begin{align*}
                \lim_{N\to\infty} \p\bigg[\sup_{t\in[0,T]}\Big|\frac{1}{N}\sum_{j\in \Z}h_j(\xN_t)f(j/N) - \int_\R H(\phi_t)(x)f(x)\dd x \Big|>\varepsilon \bigg]& \\
                \leq \lim_{N\to\infty} \p\bigg[\sup_{t\in[0,T]}\left|\frac{1}{ N}  E(\xN_t) - \int_\R \phi_t(y)\dd y\right|>\frac{\varepsilon}{ ||\psi'||_\infty} \bigg]&,
            \end{align*}
            which goes to $0$ since $\Lambda_N(\xN_t)\toP \phi_t(x)\dd x$ uniformly in $t\in[0,T]$.
            
            Now let $f\in C_c^\infty(\R)$, then there exists $M\in \N$ such that $\text{supp}(f)\subseteq [-M,M]$, hence
            \begin{align*}
                \sum_{j\in \Z}h_j(\xN_t)f(j/N) &= \sum_{j=-NM}^{NM}h_j(\xN_t) f(j/N), \\
                \int_\R H(\phi_t)(x)f(x)\dd x &= \int_{-M}^M H(\phi_t)(x)f(x)\dd x.
            \end{align*}
            Moreover, using integration by parts we have that
            \begin{align*}
                \int_{-M}^M H(\phi_t)(x)f(x)\dd x = \frac{1}{\gamma} \int_{-M}^M \psi\left(\gamma \int_x^\infty \phi_t(y)\dd y\right) f'(x)\dd x.
            \end{align*}
            Define 
            \[
                S^{(N)}= \frac{1}{N}\sum_{j=-NM}^{NM} h_j(\xN_t)f(j/N) - \frac{1}{\gamma}\int_{-M}^M \psi\left(\gamma \int_x^\infty \phi_t(y)\dd y\right) f'(x)\dd x,
            \]
            then we have to show that $\sup_{t\in [0,T]} |S^{(N)}|$ converges to $0$ in probability. For convenience of notation, in what follows we suppress the dependence on $N$ and write 
            $S= S^{(N)}$. We have $|S|\leq |S_1| + |S_2| + |S_3|$, where 
            \begin{align*}
                 S_1&=\frac{1}{N}\sum_{j=-NM}^{NM}h_j(\xN_t)f(j/N) - \frac{1}{N\gamma}\sum_{j=-NM}^{NM}\psi\Big(\frac{\gamma}{N} E_j(\xN_t)\Big)f'(j/N) , \\
                 S_2&= \frac{1}{N\gamma}\sum_{j=-NM}^{NM}\psi\Big(\frac{\gamma}{N} E_j(\xN_t)\Big)f'(j/N)   - \frac{1}{N\gamma}\sum_{j=-NM}^{NM}\psi\left(\gamma \int_{j/N}^\infty \phi_t(y)\dd y\right)f'(j/N) ,\\
                 S_3&= \frac{1}{N\gamma}\sum_{j=-NM}^{NM}\psi\left(\gamma \int_{j/N}^\infty \phi_t(y)\dd y\right)f'(j/N)  - \frac{1}{\gamma}\int_{-M}^M \psi\left(\gamma \int_x^\infty \phi_t(y)\dd y\right) f'(x)\dd x.
            \end{align*}
            We will show that the supremum over $t\in[0,T]$ of each of these terms converges to $0$ in probability. Let us start with $S_1$. 
            Using the explicit expression of $h_j^N$ and doing a summation by parts,
            \begin{align*}
                \frac{1}{N}\sum_{j=-NM}^{NM}h_j^N(\xN_t)f(j/N) = \frac{1}{N\gamma} \sum_{j=-NM}^{NM} \psi\left(\frac{\gamma}{N}E_j(\xN_t)\right)N\left[f(j/N)-f((j-1)/N)\right].
            \end{align*}
            By Taylor's theorem, we have a $\xi_j \in [(j-1)/N,j/N]$ such that  
            \[
                f(j/N)-f((j-1)/N) = \frac{f'(j/N)}{N} - \frac{f''(\xi_j)}{2N^2}.
            \]
            Therefore,
            \begin{align*}
                |S_1| \leq \frac{1}{2\gamma N^2} \sum_{j=-NM}^{NM} \left|\psi\left(\frac{\gamma}{N}E_j(\xN_t) \right)f''(\xi_j) \right| \leq  \frac{(2M+1)||\psi||_\infty ||f''||_\infty}{2\gamma N}.
            \end{align*}
            Therefore, 
            \begin{align*}
                \lim_{N\to\infty}\p\left(\sup_{t\in[0,T]}|S_1|>\varepsilon\right) &\leq \lim_{N\to\infty} \p\left(\frac{(2M+1)||\psi||_\infty ||f''||_\infty}{2\gamma  N} > \varepsilon\right) =0.
            \end{align*}
            For $S_2$ define
            \[
                A_{j,N}(\xN_t) =\sup_{t\in[0,T]} \left|\frac{1}{N} E_j(\xN_t)   -  \int_{j/N}^\infty \phi_t(y)\dd y\right|.
            \]
            Note that $A_{j,N}\leq E$ by \eqref{eq:Eass}. Let $\delta>0$, using again that $|\psi(x)-\psi(y)|\leq |x-y|\cdot ||\psi'||_\infty$, we have by Markov's inequality,
            \begin{align*}
                \p\left(\sup_{t\in[0,T]}|S_2|>\varepsilon\right) &\leq \frac{||\psi'||_\infty ||f'||_\infty}{N \varepsilon} \sum_{j=-NM}^{NM} \E\left[ A_{j,N}\left(\mathbf{1}_{A_{j,N}\leq \delta}+\mathbf{1}_{A_{j,N}>\delta}\right)\right] \\
                &\leq \frac{||\psi'||_\infty ||f'||_\infty}{N \varepsilon} \sum_{j=-NM}^{NM} \left(\delta + \E\left[ A_{j,N}\mathbf{1}_{A_{j,N}>\delta}\right]\right) \\
                &\leq \frac{||\psi'||_\infty ||f'||_\infty}{\varepsilon}  \left[(2M+1)\delta + \frac{E}{N}\sum_{j=-NM}^{NM} \p\left( A_{j,N}>\delta\right) \right].
            \end{align*}
            Since 
            \[
                \lim_{N\to\infty}\frac{E}{N}\sum_{j=-NM}^{NM} \p\left( A_{j,N}>\delta\right)=0
            \]
            by Lemma~\ref{lem:convpartialenergy} with $\mu^N_t = \Lambda_N(\xN_t)$ and $\mu_t = \phi_t(x)\dd x$, we can first let $N\to \infty$ and then $\delta \to 0$ to obtain $$\lim_{N\to\infty}\p\left(\sup_{t\in[0,T]}|S_2|>\varepsilon\right)=0.$$
            
            Note that $S_3$ is a deterministic expression where an integral over $$f_t(x)=\psi\left(\gamma \int_x^\infty \phi_t(y)\dd y\right) f'(x)$$ is approximated by a Riemann sum. Note that the $f_t$ are uniformly equicontinuous, since all $f_t$ are Lipschitz with constant $$K=||f''||_\infty ||\psi||_\infty + \sup_{t\in[0,T]}\gamma ||\phi_t||_\infty ||\psi'||_\infty||f'||_\infty.$$
            Therefore, the convergence is uniform in $t\in[0,T]$,
            \begin{align*}
                |\gamma S_3| &= \left|\sum_{j=-NM}^{NM}\left(\frac{1}{N}f_t(j/N)  - \int_{j/N}^{(j+1)/N} f_t(x)\dd x\right) \right|  \\
                &= \left|\sum_{j=-NM}^{NM}\int_{j/N}^{(j+1)/N}f_t(j/N)- f_t(x)\dd x \right| \\
                &\leq \sum_{j=-NM}^{NM}\int_{j/N}^{(j+1)/N}|f_t(j/N)- f_t(x)|\dd x  \\
                &\leq K \sum_{j=-NM}^{NM}\int_{j/N}^{(j+1)/N}|j/N- x|\dd x \\
                &\leq \frac{K}{N} \sum_{j=-NM}^{NM}\int_{j/N}^{(j+1)/N}1\dd x \\
                &\leq \frac{K(2M+1)}{N}.
            \end{align*}
            Thus for $N\geq K(2M+1)/\gamma\varepsilon$, we have 
            \[
                \p\left(\sup_{t\in[0,T]} |S_3| > \varepsilon \right) =0. \qedhere
            \]
        \end{proof}
        \begin{remark}\*
            \begin{enumerate}[label=(\roman*)]
                \item As said in Section~\ref{sec:intro},  the uniformly bounded and tight conditions on $\phi_t$ are satisfied when $t\mapsto \int_\R f(x)\phi_t(x)\dd x$ is continuous for all $f\in C_c^\infty(\R)$ and $f=1$.
                \item If we have an interacting particle system where the total energy is conserved, such as BEP, ABEP and dynamic ABEP, assumption~\eqref{eq:Eass} for $t>0$ follows from the same inequality at $t=0$.
            \end{enumerate}
        \end{remark}
    \section{Hydrodynamic limit of BEP} \label{sec:BEP}
In this section we prove the hydrodynamic limit of BEP.
We follow the classical methodology, based on the Dynkin martingale associated to the density field
(see e.g. \cite{sepa}, \cite{Demasi}). See also \cite{fran} for a recent proof of the hydrodynamic limit for a class of models including SIP and BEP on the torus. The reader acquainted with this method can skip this section which we include only for the sake of self-containedness.

\begin{theorem}[BEP]\label{thm:BEP}
Assume that the sequence of empirical measures corresponding to (random) initial configurations $\{\xN_0\}_{N\geq1}$ converges in probability to a deterministic measure which is absolutely continuous with respect to the Lebesgue measure,
    \[
        \Lambda_N(\xN_0) \toP \phi_0(x)\dd x \qquad \text{in }\mathcal{M}_+.
    \]
    Moreover, assume that the second moment of $(\xN_0)_j$ is uniformly bounded in $j$ and $N$, i.e. there exists a constant $C_\mathrm{B}$ such that
    \begin{align}
            \sup_{j\in \Z, N\in \N} \E\bigg[\big(\xN_0\big)_j^2\bigg] \leq C_\mathrm{B}, \label{eq:BEPboundquadexpectation}
        \end{align}
        and the total energy grows at most linearly in $N$, i.e. there exists a constant $E_\mathrm{B}$ such that for all $N$ we almost surely have
        \begin{align}
            \frac{1}{N}\sum_{j\in \Z} (\xN_0)_j \leq E_\mathrm{B}. \label{eq:BEPboundtotalenergy}
        \end{align}
        Then for each time $T>0$, the empirical measure of the accelerated BEP$(k)$ process started at $\xN_0$, converges in probability to the measure $\phi_t(x)\dd x$, uniformly in $t\in[0,T]$, 
        \[
        \Lambda_N(\xN_t) \toP \phi_t(x)\dd x \qquad \text{in }\mathcal{M}_+ \text{ uniformly in } t\in [0,T],
         \]
         where $$\phi_t=K_t * \phi_0$$ is the unique solution in $C([0,T],L^1(\R))\cap C^{1,2}((0,T]\times \R)$ of the heat equation 
    \begin{align*}
        	\partial_t \phi_t = k \partial_{xx}\phi_t
    \end{align*}
    satisfying 
    $$\lim_{t\downarrow 0}||\phi_t-\phi_0||_1=0.$$ 
\end{theorem}
\begin{remark} Similar to Remark~\ref{rem:muintegrable}, we have $\phi_0 \in L^1(\R) \cap L^\infty(\R)$ with $||\phi_0||_1 \leq E_\mathrm B$ and $||\phi_0||_\infty \leq \sqrt{C_\mathrm B}$. Also, since $\phi_t = K_t*\phi_0$, we have $||\phi_t||_1=||\phi_0||_1\leq E_\mathrm B$, $||\phi_t||_\infty \leq ||\phi_0||_\infty \leq  \sqrt{C_\mathrm B}$ and $\phi_t \in C^\infty(\R)$ for each $t>0$.
\end{remark}
    Proving this hydrodynamic limit is standard, follows the line of proof in \cite{sepa} and revolves around Dynkin's martingale 
    \[
            M_t^N = \langle  \Lambda_N(\xN_t),f\rangle - \langle  \Lambda_N(\xN_0),f\rangle - \int_{0}^t N^2\big[L^{\mathrm{BEP}}_{k}\langle \Lambda_N(\cdot),f\rangle\big] (\xN_s) \dd s.
        \]
   The four main ingredients are a second moment bound of $\xN_t$, vanishing of the quadratic variation of $M_t^N$, the generator $L^\mathrm{BEP}_k$ approximates the Laplacian and tightness of $\Lambda_N(\xN_t)$. To ease notation, write $f_j^N=f(j/N)$, $\nabla_{j,j+1} = \left(\partial_{x_j}-\partial_{x_{j+1}}\right)$ and $F_N(x)= \langle \Lambda_N(x),f\rangle$.
   \subsection{Second moment bound}
   \begin{proposition}\label{prop:BEPboundquadexpectationt}
        Under assumption \eqref{eq:BEPboundquadexpectation}, we have
        \begin{align*}
            \sup_{t\in [0,T],i,j\in \Z,N\in \N} \E \left[ \left(\xN_t\right)_i\left(\xN_t\right)_j \right] \leq \begin{cases}
                C_\mathrm{B} &\text{if } i\neq j, \\
                \frac{k+1}{k}C_\mathrm{B} &\text{if } i=j.
            \end{cases}
        \end{align*}
    \end{proposition} 
        \begin{proof}
            We use the duality with SIP with duality function 
            \[
                D(x,\eta)=\prod_{j\in \Z} \frac{1}{(k)_{\eta_j}} x_j^{\eta_j},
            \]
            where $(k)_n=k(k+1)\cdots(k+n-1)$ is the shifted factorial (Pochhammer symbol).
            Let $j_1(t)$ and $j_2(t)$ be the sites of the two particles at time $tN^2$ of a SIP process started from the initial condition where there are two particles at site $j$, then using Cauchy-Schwarz,
            \begin{align*}
                \frac{1}{k(k+1)}\E \left[ (\xN_t)_j^2 \right] \leq \frac{1}{k^2}\E \left[ (\xN_0)_{j_1(t)}(\xN_0)_{j_2(t)}\right]  \leq \frac{1}{k^2} \sqrt{\E \left[ (\xN_0)^{2}_{j_1(t)}\right]}\sqrt{\E \left[ (\xN_0)^{2}_{j_2(t)}\right]}.
            \end{align*}
            Therefore,
            \[
                \sup_{t\in [0,T],j\in \Z,N\in \N} \E \left[ (\xN_t)_j^2 \right] \leq \frac{k+1}{k}  \sup_{j\in \Z,N\in \N} \sqrt{\E \left[ (\xN_0)^{2}_{j}\right]} \sup_{j\in \Z,N\in \N}\sqrt{\E \left[ (\xN_0)^{2}_{j}\right]}  \leq \frac{k+1}{k} C_\mathrm{B}. 
            \]
            The case where $i\neq j$ is easier since then 
            \[
                \E \left[ (\xN_t)_i(\xN_t)_j \right] \leq \sqrt{\E \left[ (\xN_0)^{2}_{j_1(t)}\right]}\sqrt{\E \left[ (\xN_0)^{2}_{j_2(t)}\right]}. \qedhere
            \]
        \end{proof}    
    \subsection{Quadratic variation of $M_t^N$ vanishes} 
    The case with the test function $f=1$ is trivial, therefore we only have to show it for $f\in C_c^\infty(\R)$ with $\operatorname{supp}(f)\subseteq [-M,M]$. Using Doob's maximal $L^2$-inequality  we have 
    $$\E\left[\sup_{t\in[0,T]} (M_t^N)^2\right]\leq 4 \E\left[(M_T^N)^2\right].$$  
    Since $M_0^N=0$, we have $\E[(M_T^N)^2] = \E[\langle M^N\rangle_T]$. The predictable quadratic variation can be expressed using the carr\'e du champ operator,
            \[
                 \langle M^N\rangle_T  = N^2\int_0^T   \big[L^{\mathrm{BEP}}_{k}F_N^2\big](\xN_t) - 2F_N(\xN_t)\big[L^{\mathrm{BEP}}_{k}F_N \big](\xN_t)  \dd t .
            \]
    By the chain rule, we have
            \begin{align*}
                &\nabla_{j,j+1}(F_N^2) = 2 F_N \nabla_{j,j+1}F_N,\\
                &\nabla_{j,j+1}^2(F_N^2)= 2 \big(\nabla_{j,j+1}F_N\big)^2 +2F_N \nabla_{j,j+1}^2F_N.
            \end{align*}
            Therefore, using $\partial_{x_j}F_N(x) = \frac1Nf_j$,
            \begin{align}
                \begin{split}N^2(L^\mathrm{BEP}_k(F_N^2)-2F_NL^\mathrm{BEP}_k(F_N)) &= \sum_{|j|\leq NM+1} 2(\xN_t)_j(\xN_t)_{j+1}  N^2\big(\nabla_{j,j+1}F_N\big)^2 \\
                &=   \sum_{|j|\leq NM+1} 2(\xN_t)_j(\xN_t)_{j+1}  \big(f^N_j-f^N_{j+1}\big)^2.\end{split} \label{eq:cancelfirstorder}
            \end{align}
            Since $f\in C_c^\infty(\R)$, we have 
            \[
                |f_j^N-f_{j+1}^N| \leq \frac{1}{N} ||f'||_\infty.
            \]
            Therefore, using Fubini and Proposition~\ref{prop:BEPboundquadexpectationt}, we obtain
            \begin{align*}
                \E \left[  (M_T^N)^2 \right] &\leq  \frac{2}{N^2}  \sum_{|j|\leq NM+1}||f'||_\infty^2  \int_0^T\E\left[(\xN_t)_j(\xN_t)_{j+1}\right] \dd t\\
                &\leq \frac{(k+1)(4M+6)||f'||_\infty^2 C_\mathrm{B}T}{kN},
            \end{align*}
            which goes to $0$ as $N\to\infty$.
    \subsection{Generator approximates Laplacian}
    Since 
            \[
                \nabla_{j,j+1}^2F_N = 0,
            \]
            we only have to bound the first order drift term of $\BEP{k}$. Using summation by parts and a Taylor approximation, we have
    \begin{align*}
        N^2(L^\mathrm{BEP}_kF_N)(\xN_t) &=  N k\sum_{j=-NM}^{NM} ((\xN_t)_{j+1}-(\xN_t)_j)(f_j-f_{j+1}) \\
        &= N k\sum_{|j|\leq NM+1} (\xN_t)_j ( f_{j-1}^N - 2 f_j^N + f_{j+1}^N) \\
        &= N k \sum_{|j|\leq NM+1} (\xN_t)_j \left(\frac{1}{N^2}(f_j^N)'' + \frac{A}{12N^4} \right), 
    \end{align*}
    where $|A| \leq ||f^{(\text{iv})}||_\infty$. 
    Therefore, by conservation of energy, we have
    \begin{align*}
        \E\left[ \sup_{t\in[0,T]} \left|N^2(L^\mathrm{BEP}_kF_N)(\xN_t)- k\langle \Lambda_N(\xN_t),f''\rangle\right| \right] &\leq \frac{k||f^{(\text{iv})}||_\infty}{N^3} \E\left[\sup_{t\in[0,T]}  \sum_{|j|\leq NM+1} (\xN_t)_j\right] \\
        &\leq \frac{k||f^{(\text{iv})}||_\infty E_\mathrm{B}}{N^2}.
    \end{align*}
    \subsection{Tightness}
    We use Mitoma's theorem to reduce to tightness of the sequence $F_N\big(\xN_t\big)=\langle  \Lambda_N(\xN_t),f\rangle$ in $\mathcal{D}([0,T],\R)$ for each $f\in C^\infty_c(\R)$ and $f=1$. Fix such a test function. This, in turn, we will prove using Aldous' condition, i.e.  
        \begin{enumerate}[label=(\roman*)]
            \item Compact containment: for every $\varepsilon>0$, there exists a compact set $K_{\varepsilon}\subset \R $ such that
            \[
                \sup_{N\in \N}\p \big[F_N\big(\xN_t\big)\not\in K_\varepsilon \text{ for some }t\in [0,T]\big]\leq \varepsilon.
            \]
            \item Modulus of continuity, for which it is enough to show (\cite[Chapter 3]{Bi}),
            \[
                \lim_{\delta \downarrow 0} \limsup_{N\to\infty} \p\left[ \sup_{\substack{0\leq s\leq t\leq T\\ |t-s|\leq \delta}} |F_N(\xN_t)-F_N(\xN_s)| > \varepsilon \right]=0 . 
            \]
        \end{enumerate}
        Condition (i) follows directly from $||f||_\infty < \infty$, conservation of energy and assumption \eqref{eq:BEPboundtotalenergy} on the total energy. Indeed, let $K_\varepsilon=\{x\in \R: |x|\leq ||f||_\infty E_\mathrm{B}\}$, then 
        \[
            \big|F_N\big(\xN_t\big)\big| \leq \frac1N\sum_{j=-NM}^{NM} |f(j/N)|(\xN_t)_j \leq ||f||_\infty E_\mathrm{B}.
        \]
        Condition (ii) is trivial for $f=1$ by the conservation of energy, thus let us take $f\in C_c^\infty(\R)$. Then there is an $M>0$ such that $\text{supp}(f)\subseteq[-M,M]$. Using Dynkin's Martingale, we have
        \begin{align*}
            F_N(\xN_t)-F_N(\xN_s) = M_t^N - M_s^N + \int_{s}^t N^2 \big[L^{\mathrm{BEP}}_{k}F_N\big] (\xN_u) \dd u.
        \end{align*}
        Therefore, it is enough that both 
        \begin{align}
            &\lim_{\delta \downarrow 0} \limsup_{N\to\infty} \p\left[ \sup_{\substack{0\leq s\leq t\leq T\\ |t-s|\leq \delta}} |M_t^N-M_s^N| > \varepsilon \right]=0, \label{eq:BEPmodcont1}\\ 
            &\lim_{\delta \downarrow 0} \limsup_{N\to\infty} \p\left[ \sup_{\substack{0\leq s\leq t\leq T\\ |t-s|\leq \delta}} \left| \int_{s}^t N^2 \big[L^{\mathrm{BEP}}_{k}F_N\big] (\xN_u) \dd u\right| > \varepsilon \right]=0, \label{eq:BEPmodcont2}
        \end{align}
            For \eqref{eq:BEPmodcont1}, we first estimate the difference $|M_t^N-M_s^N|$ by the variation of $M_T^N$, which we already estimated. Using Markov's inequality,  $(a+b)^2\leq 2a^2 + 2b^2$ and Doob's maximal $L^2$ inequality we obtain
            \begin{align*}
                \p\left[ \sup_{\substack{0\leq s\leq t\leq T\\ |t-s|\leq \delta}} |M_t^N-M_s^N| > \varepsilon \right] &\leq \frac{1}{\varepsilon^2}\E \left[ \sup_{\substack{0\leq s\leq t\leq T\\ |t-s|\leq \delta}} (M_t^N-M_s^N)^2 \right] \\
                &\leq \frac{2}{\varepsilon^2}\E \left[ \sup_{\substack{0\leq s\leq t\leq T\\ |t-s|\leq \delta}} (M_t^N)^2+(M_s^N)^2 \right] \\
                & \leq \frac{4}{\varepsilon^2}\E \left[ \sup_{0\leq t\leq T} (M_t^N)^2 \right]\\
                &\leq \frac{16}{\varepsilon^2}\E \left[  (M_T^N)^2 \right],
            \end{align*}
            which goes to $0$ as $N\to \infty$.

            \medskip
            
            Let us now show \eqref{eq:BEPmodcont2}. Using the same computation as for the generator approximation, we have
            \begin{align*}
                \left| N^2\big[L^\mathrm{BEP}_k F_N \big](\xN_u)\right| \leq k||f''||_\infty E_\mathrm{B}.
            \end{align*}
            The limit \eqref{eq:BEPmodcont2} now follows from letting $\delta \downarrow 0$. 
    \section{Proof of hydrodynamic limit of ABEP (Theorem~\ref{thm:ABEP})} \label{sec:ABEP}
    In this section we exploit the general micro-to-macro transformation result from Section~\ref{sec:microtomacro} to prove the hydrodynamic limit of ABEP. As explained in Remark~\ref{rem:muintegrable}(i), we only have to prove the case $\gamma>0$. We proceed in five steps.
    \begin{enumerate}
        \item Since by assumption the empirical measure $\Lambda_N(\yN_0)\toP \mu_0(y)\dd y$, we can use Theorem~\ref{thm:microtomacro} with $T=0$ to show that the microscopic transformation $g$ converges to the macroscopic transformation $G$ (see Table~\ref{tab:micro-macro-transformations}), $$\Lambda_N(g(\yN_0))\toP G(\mu_0)(x)\dd x.$$
        \item Since $\xN_t=g(\yN_t)$ follows BEP dynamics, we can use Theorem~\ref{thm:BEP} to show that $\Lambda_N\left(g(\yN_t)\right)$ converges to the solution $\phi_t(x)$ of the heat equation with initial condition $\phi_0(x)=G(\mu_0)(x)$, uniformly in $t\in[0,T]$. 
        \item We do step (1) the other way around. We use Theorem~\ref{thm:microtomacro} to show that the empirical measure of $\yN_t= g^{-1}(\xN_t)$ converges to $\mu_t(y)=G^{-1}(\phi_t)(y)$, uniformly in $t\in [0,T]$. Here, $G^{-1}$ (see Table~\ref{tab:micro-macro-transformations}) is the macroscopic version of the transformation $g^{-1}$.
        \item We show that $G(\mu_t)(x)=G(G^{-1}(\phi_t))(x)=\phi_t(x)$. Therefore, $G(\mu_t)(x)$ solves the heat equation with initial condition $G(\mu_0)(x)$.
        \item From step 4, we show that $\mu_t(x)$ solves the viscous Burgers equation with initial condition $\mu_0(x)$.
    \end{enumerate}
        For step (1), we take $T=0$ and $\psi(a)=-e^{-a}$ in Theorem~\ref{thm:microtomacro}. The latter gives 
        \begin{align*}
            h_j^N(y) &= g_j(y), \\
            H(\mu)(x) &= G(\mu)(x).
        \end{align*}
        Therefore, we obtain the following corollary, using assumptions \eqref{eq:ABEPboundquadexpectation} and \eqref{eq:ABEPboundtotalenergy} and the convergence $\Lambda_N(\yN_0)\toP \mu_0(y)\dd y$ of initial configurations.
        \begin{corollary}
             The empirical measure of the transformation $g$ of the initial configurations $\yN_0$ converges to the deterministic measure $G(\mu_0)(x)\dd x$, $$\Lambda_N(g(\yN_0)) \toP G(\mu_0)(x)\dd x.$$
        \end{corollary}
        Since $\xN_t=g(\yN_t)$ follows BEP dynamics, we will use Theorem~\ref{thm:BEP} to obtain that the empirical measure of $g(\yN_t)$ converges to the solution of the heat equation, completing step (2). We first need to check the bounds on the second moment and total energy of $g(\yN_0)$, which follows from the explicit form of $g$ and assumptions \eqref{eq:ABEPboundquadexpectation}, \eqref{eq:ABEPboundtotalenergy} on the initial states $\yN_0$.
        \begin{proposition} \label{prop:ABEPbounds}
            Assume that \eqref{eq:ABEPboundquadexpectation} and \eqref{eq:ABEPboundtotalenergy} hold, then we have 
            \begin{align*}
                \sup_{j\in\Z,N\in\N} \E\big[g_j(\yN_0)^2\big] &\leq C_\mathrm{A}, \\
                \frac{1}{N} \sum_{j\in \Z}g_j(\yN_t) &\leq \frac{1-e^{-\gamma E_\mathrm{A}}}{\gamma}
                \ \text{ for all}\ t\geq 0.
            \end{align*}
        \end{proposition}
        \begin{proof}
           By the mean value theorem $|e^{-x}-e^{-y}|\leq |x-y|$ for $x,y\geq 0$, hence
            \[
                g_j(\yN_0) \leq (\yN_0)_j.
            \]
            The required inequality now follows from  Assumption \eqref{eq:ABEPboundquadexpectation}. For the second inequality, note that by telescoping we have for all $n_1 < n_2$,
            \[
                \frac{1}{N}\sum_{j=n_1}^{n_2}g_j(\yN_0) = \frac{1}{\gamma} \big( e^{-\frac{\gamma}{N} E_{n_2+1}(\yN_0)} - e^{-\frac{\gamma}{N} E_{n_1}(\yN_0)}\big) \leq \frac{1-e^{-\gamma E_\mathrm{A}}}{\gamma}.
            \]
            Since the total energy is conserved, the proof is finished.
        \end{proof}
        \begin{corollary}
            For all $T>0$, 
            $$\Lambda_N(g(\yN_t))\toP \phi_t(x)\dd x \quad \text{in } \mathcal{M}_+ \text{ uniformly in } t\in[0,T],$$
            where $\phi_t$ is the unique solution in $C([0,T],L^1(\R))\cap C^{1,2}((0,T]\times \R)$ of the heat equation 
            $$\partial_t \phi_t = k \partial_{xx}\phi_t, $$
           satisfying  $$\lim_{t\downarrow 0}||\phi_t-G(\mu_0)||_1=0.$$
           Moreover, we have 
            \begin{align*}
                ||\phi_t||_1 &= ||\phi_0||_1 \leq \frac{1-e^{-\gamma E_\mathrm{A}}}{\gamma},\\
                ||\phi_t||_\infty &\leq ||\phi_0||_\infty \leq \sqrt{C_\mathrm A}.
            \end{align*}
        \end{corollary}
        Now we would like to show that the empirical measure of $\yN_t$ converges to an absolutely continuous measure $\mu_t(y)\dd y$, which would satisfy $G(\mu_t)(x)=\phi_t(x)$ by Theorem~\ref{thm:microtomacro} and uniqueness of limits. This follows from Theorem~\ref{thm:microtomacro} where we go in the other direction. That is, we use it with $\xN_t = g(\yN_t)$ and $\psi(a)= - \log(1-a)$, which is well-defined by the second inequality of Proposition~\ref{prop:ABEPbounds}. This gives
        \begin{align*}
            h_j^N(x) = g^{-1}_j(x), \\
            H(\phi)(y) := G^{-1}(\phi)(y). 
        \end{align*}
        
        Note that the assumptions \eqref{eq:ubddtight} and \eqref{eq:Eass} are satisfied. For the first, it is enough to show that $t\mapsto \int_{\R} f(x)\phi_t(x)\dd x$ is continuous for all $f\in C_c^\infty(\R)$ and $f=1$. This follows from the fact that $t\mapsto \phi_t$ is continuous in $L^1(\R)$ via
        \[
            \left|\int_\R (\phi_t(x)-\phi_s(x))f(x)\dd x \right| \leq  ||\phi_t-\phi_s||_{1} ||f||_\infty. 
        \]
        For the second, the energy bound on $\xN_t=g(\yN_t)$ is satisfied by Proposition~\ref{prop:ABEPbounds}. Thus, the empirical measure of $g^{-1}(\xN_t)=g^{-1}(g(\yN_t))=\yN_t$ converges to $G^{-1}(\phi_t)(y)\dd y$, uniformly in $t\in[0,T]$, finishing step (3).
        \begin{corollary}
            Define $\mu_t(y) =G^{-1}(\phi_t)(y)$, then 
            $$ \Lambda_N(\yN_t) \toP \mu_t(y) \dd y, \quad \text{in } \mathcal{M}_+\text{ uniformly in }t\in[0,T].$$
        \end{corollary}
        For step (4), we show that $G(\mu_t)(x) = G(G^{-1}(\phi_t))(x)= \phi_t(x)$. Recall that both $G$ and $G^{-1}$ are of the form
        \[
            H(\phi_t)(x)= -\frac{1}{\gamma} \partial_x \psi\left(\gamma \int_x^\infty \phi_t(y)\dd y\right).
        \]
        Let $\psi_{\log}(a)=-\log(1-a)$ and $\psi_e(a)=-e^{-a}$, then we have
        \begin{align*}
            G(G^{-1}(\phi_t))(x) &= -\frac{1}{\gamma} \partial_x\psi_e\left( \gamma\int_x^\infty G^{-1}(\phi_t)(y)\dd y\right) \\
            &= -\frac{1}{\gamma} \partial_x\psi_e\left( -\int_x^\infty \partial_y\psi_{\log}\left(\gamma  \int_y^\infty \phi_t(z)\dd z\right) \dd y\right) \\
            &= -\frac{1}{\gamma} \partial_x\psi_e\left( \psi_{\log}\left(\gamma \int_x^\infty \phi_t(z)\dd z\right) - \psi_{\log}(0)\right) \\
            &= \frac{1}{\gamma} \partial_x \left(1-\gamma\int_x^\infty \phi_t(z)\dd z\right) \\
            &= \phi_t(x),
        \end{align*}
        where we used the Lebesgue differentiation theorem in the last step. Therefore, $G(\mu_t)(x)=\phi_t(x)$ satisfies the heat equation. From this, we can show that $\mu_t$ satisfies the viscous Burgers equation. We need the following two lemmas. The first shows that $$\int_x^\infty G(\mu_t)(y) \dd y = \frac{1-e^{-\gamma \int_x^\infty \mu_t(y)\dd y}}{\gamma}$$ also satisfies the heat equation with initial condition $$\int_x^\infty G(\mu_0)(y)\dd y = \frac{1-e^{-\gamma\int_x^\infty \mu_0(y)\dd y}}{\gamma}.$$ By linearity of the heat equation, $$J_t(x)=e^{-\gamma \int_x^\infty \mu_t(y)\dd y}$$ then also satisfies the heat equation with initial condition $$J_0(x)=e^{-\gamma \int_x^\infty \mu_0(y)\dd y}.$$ The second lemma shows that then $\mu_t(y)$ satisfies viscous Burgers equation.
        \begin{lemma}\label{lem:derheateqn}
            Let $k>0$ and $\phi_0 \in L^1(\R)$. Suppose $\phi_t(x)$ is the unique solution in $C([0,T],L^1(\R))\cap C^{1,2}((0,T]\times \R)$ of the heat equation 
            \begin{align*}
                &\partial_{t}\phi_t(x) = k \partial_{xx}\phi_t(x), \qquad (t,x)\in (0,T]\times \R,
            \end{align*}
            satisfying 
            $$\lim_{t\downarrow 0} ||\phi_t - \phi_0||_1=0 .$$
            Then the function $$\Phi_t(x)=\int_x^\infty \phi_t(y)\dd y$$ is the unique solution in $C([0,T],\BUCR)\cap C^{1,2}((0,T]\times \R)$ of the heat equation satisfying
            $$\lim_{t\downarrow 0} ||\Phi_t - \Phi_0||_\infty=0,$$
            where $\Phi_0(x)=\int_x^\infty \phi_0(y)\dd y$.
        \end{lemma}
        \begin{proof} 
        We will show that $\Phi_t = K_t*\Phi_0$, the result then follows from uniqueness and strong continuity of the heat semigroup on $\BUCR$. \\
        
        \noindent Since $\phi_t$ solves the heat equation, we have $$\phi_t=K_t*\phi_0,$$
        where $K_t$ is the heat kernel. Using the definition of $\Phi_t$ combined with Fubini's theorem, 
        \begin{align*} 
            \Phi_t(x) &=\int_x^\infty\int_{\mathbb{R}} K_t(y-z)\phi_0(z)\,\mathrm{d}z\,\mathrm{d}y\\ 
            &=\int_{\mathbb{R}}\phi_0(z) \left(\int_x^\infty K_t(y-z)\,\mathrm{d}y\right)\mathrm{d}z. 
            \end{align*} 
            The interchange of integrals is justified since by Young's inequality \[ \int_{\mathbb{R}}|\phi_0(z)| \left(\int_x^\infty K_t(y-z)\,\mathrm{d}y\right)\mathrm{d}z \leq \|\phi_0\|_{L^1}. \] 
            On the other hand, we have
            \begin{align*} 
            K_t*\Phi_0&=\int_{\mathbb{R}}K_t(x-w) \left(\int_w^\infty\phi_0(z)\,\mathrm{d}z\right)\mathrm{d}w\\ 
            &=\int_{\mathbb{R}}\phi_0(z) \left(\int_{-\infty}^zK_t(x-w)\,\mathrm{d}w\right)\mathrm{d}z \\
            &= \int_{\mathbb{R}}\phi_0(z) \left(\int_x^\infty K_t(y-z)\,\mathrm{d}y\right) \dd z,
            \end{align*} 
            where we used the change of variables $y=x-w+z$ in the last step. It follows that $\Phi_t=K_t*\Phi_0$. 
        \end{proof}
        \begin{lemma}
            Let $\gamma\neq 0$, $\mu_0\in L^1(\mathbb{R})$ and define \[ J_0(x) = \exp\left(-\gamma\int_x^\infty \mu_0(y)\,\dd y\right). \] Suppose that $J_t$ is the unique solution in $C([0,T],\BUCR)\cap C^{1,2}((0,T]\times \R)$ of the heat equation 
            \[ 
                \partial_t J_t = k\partial_{xx}J_t, \qquad k>0, 
            \] 
            satisfying 
            \[
                \lim_{t\downarrow 0} || J_t - J_0||_\infty = 0.
            \]
            Then \[ \mu_t(y):=\frac{1}{\gamma}\partial_y\log J_t(y) \] is the unique solution in $C([0,T],L^1(\R))\cap C^{1,2}((0,T]\times \R)$ of the viscous Burgers equation 
           \begin{align*}
                \partial_t\mu_t = k\partial_{yy}\mu_t + k \gamma \partial_y(\mu_t^2),\qquad (t,y)\in (0,T]\times \R,
            \end{align*}
            satisfying 
            \begin{align*}
                \lim_{t\downarrow 0} ||\mu_t - \mu_0||_1 =0.
            \end{align*}
        \end{lemma} 
        Uniqueness follows from a similar argument as in the proof of Theorem~\ref{thm:uniquegenbur}.
        \begin{proof} Since 
        \[
            J_t = K_t * J_0,
        \]
        we have 
        \[
            J_t(x) \geq \inf_{x\in \R} J_0(x) \geq e^{-|\gamma| ||\mu_0||_1}. 
        \]
        Thus, $u_t(x)=\log J_t(x)$ is well defined. For $t>0$, we use $J_t=e^{u_t}$ to obtain
        \[ 
            \partial_t J_t=J_t\,\partial_tu_t 
        \] 
        and 
        \begin{align*}
            \partial_xJ_t &= J_t \partial_xu_t, \\
            \partial_{xx}J_t &= J_t\left(\partial_{xx}u_t+\Big(\partial_xu_t\Big)^2\right). 
        \end{align*}
        Substitution into the heat equation and division by $J_t>0$ give 
        \[ 
            \partial_tu_t = k\left(\partial_{xx}u_t+(\partial_xu_t)^2\right). 
        \] 
        Differentiating this identity with respect to $x$ yields 
        \[ 
            \partial_t(\partial_xu_t) = k\,\partial_{xx}(\partial_xu_t) + k\,\partial_x\left((\partial_xu_t)^2\right), 
        \] 
        where we could interchange the partial derivatives with respect to $t$ and $x$ since $u_t$ is smooth in both variables for $t>0$. Because 
        \[ 
            \partial_xu_t=\gamma\mu_t, 
        \] 
        we have 
        \[ 
            \gamma\,\partial_t\mu_t = k\gamma\,\partial_{xx}\mu_t + k\gamma^2\,\partial_x(\mu_t^2).
        \] 
        Now divide by $\gamma\neq0$. \\
        To prove convergence of the initial condition, let $t>0$, then 
        \[
            \mu_t = \frac{ \partial_y J_t }{\gamma J_t} = \frac{ \partial_y(K_t*J_0) }{\gamma J_t} = \frac{ K_t * \partial_y J_0 }{\gamma J_t},
        \]
        where 
        \[
            \partial_y J_0 = \gamma \mu_0 J_0
        \]
        is well-defined by the Lebesgue differentiation theorem. Therefore,
        \begin{align*}
            \left| \mu_t(x) - \mu_0(x) \right| &\leq \frac{\int_\R K_t(x-y)J_0(y)|\mu_0(y)-\mu_0(x)|\dd y}{J_t(x)} \\
            &\leq e^{2|\gamma|||\mu_0||_1} \int_\R K_t(x-y)|\mu_0(y)-\mu_0(x)|\dd y,
        \end{align*}
        where we used that 
        \[
            e^{-|\gamma| ||\mu_0||_1} \leq J_t \leq  e^{|\gamma| ||\mu_0||_1} .
        \]
        The claim
        \[
            \lim_{t\downarrow 0} ||\mu_t - \mu_0||_1 =0
        \]
         now follows from the same estimate required for strong continuity of the heat semigroup on $L^1(\R)$. 
        \end{proof}
    \section{Proof of hydrodynamic limit of dynamic ABEP (Theorem~\ref{thm:dynABEP})} \label{sec:dynABEP}
       Similar to the proof in the previous section, we will proceed in five steps.
    \begin{enumerate}
        \item Since by assumption the empirical measure $\Lambda_N(\zN_0)\toP \nu_0(z)\dd z$, we can use Theorem~\ref{thm:microtomacro} with $T=0$ to show that the microscopic transformation $\gd$ converges to the macroscopic transformation $\Gd$ (see Table~\ref{tab:micro-macro-transformations}), $$\Lambda_N(\gd(\zN_0))\toP \Gd(\nu_0)(x)\dd x.$$
        \item Since $\xN_t=\gd(\zN_t)$ follows BEP dynamics, we can use Theorem~\ref{thm:BEP} to show that $\Lambda_N\left(\gd(\zN_t)\right)$ converges to the solution $\phi_t(x)$ of the heat equation with initial condition $\phi_0(x)=\Gd(\nu_0)(x)$, uniformly in $t\in[0,T]$.
        \item We do step (1) the other way around. We use Theorem~\ref{thm:microtomacro} to show that the empirical measure of $(\gd)^{-1}(\xN_t)=(\gd)^{-1}(\gd(\zN_t))=\zN_t$ converges to $\nu_t(z)=(\Gd)^{-1}(\phi_t)(z)$ uniformly in $t\in[0,T]$, where $(\Gd)^{-1}$ is the macroscopic version of the transformation $(\gd)^{-1}$, see Table~\ref{tab:micro-macro-transformations}.
        \item We show that $\Gd(\nu_t)(x)=\Gd((\Gd)^{-1}(\phi_t))(x)=\phi_t(x)$. Therefore, $\Gd(\nu_t)(x)$ solves the heat equation. 
        \item From this, we show that $V_t(z)=\int_z^\infty \nu_t(x)\dd x$ solves the non-compact dynamic deterministic KPZ equation.
    \end{enumerate}
        For step (1), we use Theorem~\ref{thm:microtomacro} with $\psi(a)=\cosh(\rho+a)$.
        \begin{corollary}\label{cor:gdynmicrotomacro}
             The empirical measure of the transformation $\gd$ of the initial configurations $\zN_0$ converges to $\Gd(\nu_0)(x)\dd x$, 
             \[\Lambda_N(\gd(\zN_0)) \toP \Gd(\nu_0)(x)\dd x.\] 
        \end{corollary}
         For step (2), we use Theorem~\ref{thm:BEP} to show that the empirical measure of $\xN_t=\gd(\zN_t)$ converges to the solution $\phi_t(x)$ of the heat equation. Again, the assumptions \eqref{eq:BEPboundquadexpectation} and \eqref{eq:BEPboundtotalenergy} follow from the explicit form of $\gd$ and assumptions \eqref{eq:dynABEPboundquadexpectation} and \eqref{eq:dynABEPboundtotalenergy}.
        \begin{proposition}\label{prop:dynABEPbounds}
            Assume that \eqref{eq:dynABEPboundquadexpectation} and \eqref{eq:dynABEPboundtotalenergy} hold, then we have
            \begin{align*}
                \sup_{j\in\Z,N\in\N} \E\big[\gd_j\big(\zN_0\big)^2\big] &\leq  C_\mathrm{D} \sinh(\rho+\gamma E_\mathrm{D})^2, \\
                \frac{1}{N} \sum_{j\in \Z}\gd_j(\zN_t) &\leq \frac{\cosh(\rho+\gamma E_\mathrm{D}) - \cosh(\rho)}{\gamma}.
            \end{align*}
        \end{proposition}
        \begin{proof}
            Let $x,y\in [\rho,\rho+\gamma E_\mathrm{D}]$ such that $x\leq y$. By the mean value theorem, there exists a $c\in (x,y)$ such that  $\cosh(x)-\cosh(y)= (x-y)\sinh(c)$. Applying this to $\gd_j$ gives
            \[
                |\gd_j(\zN_0)| = (\zN_0)_j\sinh(c) \leq (\zN_0)_j \sinh(\rho+\gamma E_\mathrm{D}).
            \]
            The first inequality now follows from Assumption \eqref{eq:dynABEPboundquadexpectation}. \\
            For the second inequality, note that by telescoping and conservation of energy we have for all $n_1 < n_2$, 
            \[
                \frac{1}{N}\sum_{j=n_1}^{n_2}\gd_j(\zN_t) \leq \frac{1}{\gamma} \big( \cosh(\rho+\gamma E_\mathrm{D}) - \cosh(\rho)\big). \qedhere
            \]
        \end{proof}
        \begin{corollary}
            For each time $T>0$, 
            $$\Lambda_N(\gd(\zN_t))\toP \phi_t(x)\dd x \quad \text{in } \mathcal{M}_+ \text{ uniformly in } t\in[0,T],$$
            where $\phi_t$ is the solution of the heat equation 
            $$\partial_t \phi_t = k \partial_{xx}\phi_t, $$
            with initial condition $\phi_0(x)=\Gd(\nu_0)(x)$. Moreover, we have 
            \begin{align*}
                ||\phi_t||_1 &= ||\phi_0||_1 \leq \frac{\cosh(\rho +\gamma E_\mathrm D) - \cosh(\rho)}{\gamma},\\
                ||\phi_t||_\infty &\leq ||\phi_0||_\infty \leq \sqrt{C_\mathrm D}\sinh(\rho+\gamma E_\mathrm D),
            \end{align*}
            and $\phi_t\in C([0,T],L^1(\R))$ and $\phi_t \in C^\infty(\R)$ for $t>0$.
        \end{corollary} 
        Then, for step (3) we use Theorem~\ref{thm:microtomacro} with $\psi(a)=\arccosh(\cosh(\rho)+a)$, which is smooth since $\rho>0$. The fact that the assumptions \eqref{eq:ubddtight} and \eqref{eq:Eass} are satisfied follows by a similar argument as in Section~\ref{sec:ABEP}.
        \begin{corollary}
            Define $\nu_t(z)= (\Gd)^{-1}(\phi_t)(z)$, then 
            \[
                    \Lambda_N(\zN_t)\toP \nu_t(z) \dd z\quad \text{in } \mathcal{M}_+ \text{ uniformly in } t\in [0,T].
            \]
        \end{corollary}
        For step (4), we show that $\Gd(\nu_t)(x) = \Gd((\Gd)^{-1}(\phi_t))(x)= \phi_t(x)$. We proceed similarly as for ABEP, but now we define
        \begin{align*}
            \psi_{\cosh}(a)&= \cosh(\rho+a),\\
            \psi_{\arccosh}(a)&=\arccosh(\cosh(\rho)+a).
        \end{align*}
        Then we have
        \begin{align*}
            \Gd((\Gd)^{-1}(\phi_t))(x) &= -\frac{1}{\gamma} \partial_x\psi_{\cosh}\left( \psi_{\arccosh}\left(\gamma \int_x^\infty \phi_t(z)\dd z\right) - \psi_{\arccosh}(0)\right) \\
            &= -\frac{1}{\gamma} \partial_x \left(\cosh(\rho)+\gamma\int_x^\infty \phi_t(z)\dd z\right) \\
            &= \phi_t(x).
        \end{align*}
        Therefore, we know that $$\Gd(\nu_t)(x)= -\frac{1}{\gamma} \partial_x\cosh\left(\rho+\gamma \int_x^\infty \nu_t(z)\dd z\right)$$ satisfies the heat equation. Using Lemma~\ref{lem:derheateqn}, we have that $$J_t(x)=\cosh\left(\rho+\gamma \int_x^\infty \nu_t(z)\dd z\right)$$ satisfies the heat equation as well, with initial condition 
        \[
            J_0(x) =  \cosh\left(\rho + \gamma \int_x^\infty \nu_0(z)\dd z\right).
        \]
        From this, we can deduce a PIDE for $\nu_t(z)$, but not a PDE. Instead, we derive a PDE for the tail potential $$V_t(z)= \int_z^\infty \nu_t(x)\dd x.$$ Essentially, this comes down to taking a hyperbolic Cole--Hopf transformation $J_t=\cosh(\rho + \gamma V_t)$, where $J_t$ satisfies the heat equation.
        \begin{lemma}\label{lem:genburg}
        Let $T,\gamma,\rho,k>0$ and $V_0\in \BUCR$ nonnegative. Suppose that $J_t^{\mathrm{dyn}}$ is the unique solution in $C([0,T],\BUCR)\,\cap\, C^{1,2}((0,T]\times \R)$ of the heat equation  
            \begin{align*}
                &\partial_{t}J_t^{\mathrm{dyn}}(x) = k \partial_{xx}J_t^{\mathrm{dyn}}(x), \qquad (t,x)\in (0,T]\times \R,
            \end{align*}
            satisfying $$\lim_{t\downarrow 0} ||J_t^{\mathrm{dyn}}-\cosh(\rho + \gamma V_0)||_\infty=0.$$
            Then $$V_t=\frac{\arccosh(J_t^{\mathrm{dyn}})-\rho}{\gamma}$$ is the unique solution in $C([0,T],\BUCR)\,\cap\, C^{1,2}((0,T]\times \R)$ of the non-compact dynamic deterministic KPZ equation,
            \begin{align*}
                \partial_tV_t(z) &= k\partial_{zz}V_t(z) + k \gamma \coth\big(\rho+\gamma V_t(z)\big)\big[\partial_zV_t(z)\big]^2, \qquad (t,z)\in (0,T]\times \R,
            \end{align*}
            satisfying $$\lim_{t\downarrow 0} ||V_t-V_0||_\infty =0 .$$
        \end{lemma}
        \begin{proof}
        Since 
        \[
            J^\mathrm{dyn}_t(x) = K_t * J^\mathrm{dyn}_0(x) \geq \inf_{x\in\R} J^\mathrm{dyn}_0(x) \geq \cosh(\rho),
        \]
        $V_t$ is well-defined, nonnegative and $J^\mathrm{dyn}_t(x)=\cosh(\rho+\gamma V_t(x))$. For $t>0$, we use the chain rule to obtain
            \begin{align*}
                \partial_tJ^{\mathrm{dyn}}_t(x) &= \gamma\sinh(\rho+\gamma V_t(x))\partial_tV_t(x),\\
                \partial_xJ^{\mathrm{dyn}}_t(x) &=\gamma  \sinh(\rho+\gamma V_t(x))\partial_xV_t(x),\\
                \partial_{xx}J^{\mathrm{dyn}}_t(x) &= \gamma  \sinh(\rho+\gamma V_t(x))\partial_{xx}V_t(x) + \gamma^2 \cosh(\rho+\gamma V_t(x))\Big[\partial_xV_t(x)\Big]^2.
            \end{align*}
            Since $J^{\mathrm{dyn}}_t$ satisfies the heat equation, we have
            \begin{align*}
                \gamma\sinh(\rho+\gamma V_t(x))\partial_tV_t(x) = k \gamma  \sinh(\rho+\gamma V_t(x))\partial_{xx}V_t(x) + k\gamma^2 \cosh(\rho+\gamma V_t(x))\Big[\partial_xV_t(x)\Big]^2.
            \end{align*}
            Now divide both sides by $\gamma\sinh(\rho+\gamma V_t(x))>0$. Existence, uniqueness and regularity of the solution now follow from Theorem~\ref{thm:uniquegenbur}. 
        \end{proof}
        \subsection{A few additional considerations}
In this subsection, we sketch a few potential directions of further research regarding the hydrodynamic limit of dynamic interacting particle systems.
\begin{enumerate}[label=(\roman*)]
    \item \textbf{Dynamic BEP}. If one takes $\sigma\to 0$ in dynamic ABEP, one obtains dynamic BEP \cite{Wa}, which has generator
    \begin{align*}
        [L^{\mathrm{BEP}_\mathsf{R}}_{k,\rho}f](z) = \sum_{j\in \Z}  c_{j}^+ c_{j+1}^-  \big[(\nabla_{j,j+1})^2f\big](z) 
        + \big( k[c_{j+1}^- - c_{j}^+] + c_{j}^+c_{j+1}^- d_j\big)\big[\nabla_{j,j+1} f\big](z) ,
    \end{align*}
    where
    \begin{align*}
        c_{j}^+ &= z_{j}\, \frac{\rho+2E_{j+1}(z)+z_{j}}{\rho+2E_{j+1}(z)}, \quad c_{j+1}^- = z_{j+1}\, \frac{\rho+2E_{j+1}(z)-z_{j+1}}{\rho+2E_{j+1}(z)}, \quad      d_j = \frac{2}{\rho+2E_{j+1}(z)}.
    \end{align*}
    Performing the same limit operation to $\gd/2\sigma$, gives that $$\gb_j(z)=\tfrac12z_j\big(\rho+ E_j(z)+E_{j+1}(z)\big)$$ is a non-local transformation connecting BEP and dynamic BEP. Under the scaling $\rho_N=2N\rho$, if 
    $$ \frac1N \sum_{j\in \Z} w_j f(j/N) \to \int_\R f(w)\omega(w)\dd w$$
    when $N\to \infty$, then formally
    \begin{align*}
        \frac{1}{N^2}\sum_{j\in \Z}\gb_j(w)f(j/N) \to& \int_\R f(w)\omega(w) \bigg(\rho+\int_w^\infty\omega (x)\dd x\bigg) \dd w\\
        &= -\frac12 \int_\R f(w) \partial_w \big(\rho+F(w)\big)^2\dd w,
    \end{align*}
    where $$F(w)=\int_w^\infty \omega(x)\dd x.$$ Consequently, $\big(\rho+F(w)\big)^2$ solves the heat equation. Therefore,
    \begin{align*}
        (\rho+F) \partial_t F &= k\partial_w \big( (\rho+F)\partial_w F\big) \\
        &= k(\rho+F)\partial_{ww}F + k\big(\partial_w F\big)^2.
    \end{align*}
    Dividing by $\rho+F >0$ then gives the PDE
    \[
        \partial_t F = k\partial_{ww}F + k\frac{\big(\partial_w F\big)^2}{\rho +F},
    \]
    which is thus linearized by the transformation $F\mapsto \big(\rho + F\big)^2$.
    It would be interesting to make this hydrodynamic limit rigorous.\\
    
    \item \textbf{Dynamic ASEP}. As explained before in the introduction, the nomenclature ``the non-compact dynamic deterministic KPZ equation'' refers to the fact that this equation arises from systems with the Lie algebra $\mathfrak{su}(1,1)$, coming from the non-compact group $SU(1,1)$. Naively, this would suggest the compact counterpart,
    \[
        \partial_tV_t = k\,\partial_{xx}V_t + k \gamma \tanh\bigl(\rho+\gamma V_t\bigr) \bigl(\partial_xV_t\bigr)^2,
    \]
    where $\rho \in \R$. This PDE is linearized by the hyperbolic-sine transformation  
    \[ 
        V\longmapsto\sinh(\rho+\gamma V). 
    \]  
    Since $\sinh$ is a bijection of $\R$, no sign
    restriction on $\rho$ is needed, in contrast with the non-compact case. 
    
    There is a clear interpolation mechanism for the PDE. Since $-1 \leq \tanh(x)\leq 1$, the coefficient of the quadratic-gradient term interpolates between the two limiting values $-k\gamma$ and $k\gamma$. Accordingly, the equation connects two deterministic KPZ-type equations with opposite signs of the nonlinearity. 
    
    For an exclusion process of spin $j\in\tfrac12\N$, i.e.\ with at most $2j$
    particles per site, the natural variable is not the tail potential
    $W_t(x)=\int_x^\infty\nu_t(y)\dd y$ but the height function
    $V_t(x)=W_t(x)+jx$, which has derivative $\partial_xV_t=j-\nu_t$. The constant profiles $\nu=0$ and $\nu=2j$ should be stationary, therefore we suspect that the correct compact PDE is given by
    \begin{align}
        &\partial_tV_t = k\,\partial_{xx}V_t + k \gamma \tanh\bigl(\rho+\gamma V_t\bigr)\Big[ \bigl(\partial_xV_t\bigr)^2-j^2\Big],\nonumber \\
        \Leftrightarrow\quad  &\partial_tW_t = k\,\partial_{xx}W_t + k \gamma \tanh\bigl(\rho+\gamma W_t+\gamma jx\bigr)\partial_xW_t \bigl(\partial_xW_t+2j\bigr),\label{eq:corCDDKPZ}
    \end{align}
    which is linearized by 
    \[
        W\mapsto e^{k \gamma^2 j^2 t}\sinh(\rho+\gamma W + \gamma jx).
    \]
    In the limit $\rho\to\pm\infty$, the density
    satisfies the viscous Burgers equation with exclusion-type flux,
    \[
        \partial_t\nu_t = k\,\partial_{xx}\nu_t
        \pm k\gamma\,\partial_x\bigl[\nu_t(2j-\nu_t)\bigr].
    \]  
    It is natural to ask whether the compact dynamic deterministic KPZ equation \eqref{eq:corCDDKPZ} arises as the hydrodynamic equation under a suitable weak-asymmetry regime of a dynamic exclusion process such as dynamic ASEP \cite{BoCo,GrWaDyn}. The latter has an algebraic interpretation in the quantum algebra $\mathcal{U}_q(\mathfrak{su}(2))$, which (formally) converges to $\mathcal{U}(\mathfrak{su}(2))$ when letting $q\to 1$, where $\mathfrak{su}(2)$ is the Lie algebra of the compact group $SU(2)$.
    \end{enumerate}
        \section{Existence and uniqueness of solutions to the non-compact dynamic deterministic KPZ equation} \label{sec:uniexi}
        \begin{theorem} \label{thm:uniquegenbur}
        Let $T,k,\gamma,\rho>0$, and $V_0 \in \operatorname{BUC}(\mathbb{R})$ with $V_0\geq 0$. Then there exists a unique solution 
        \[ 
            V \in C\bigl([0,T];\operatorname{BUC}(\mathbb{R})\bigr) \cap C^{1,2}\bigl((0,T]\times\mathbb{R}\bigr) 
        \] 
        of 
        \[ 
            \partial_t V_t(x) = k \partial_{xx}V_t(x) + k\gamma \coth\bigl(\rho+\gamma V_t(x)\bigr) \bigl[\partial_xV_t(x)\bigr]^2, \qquad (t,x)\in(0,T]\times\mathbb{R}, 
        \] 
        satisfying  
        \[ 
            \lim_{t\downarrow0} \lVert V_t-V_0\rVert_{\infty}=0. 
        \] 
        The solution is given by 
        \[ 
            V_t = \frac{1}{\gamma} \left[ \arccosh\left( K_t * \cosh\bigl(\rho+\gamma V_0\bigr)\right) -\rho \right] \geq 0, 
        \] 
        where 
        \[ 
            K_t(x) = \frac{1}{\sqrt{4\pi kt}} \exp\left(-\frac{x^2}{4kt}\right) 
        \] 
        is the heat kernel associated with $\partial_t\phi_t=k\partial_{xx}\phi_t$. \end{theorem}
        \begin{proof} Define 
        \[ 
            J_0:=\cosh(\rho+\gamma V_0) \quad\text{and}\quad J_t:=K_t*J_0. 
        \] 
        Since \(V_0\in\operatorname{BUC}(\mathbb{R})\), we have \(J_0\in\operatorname{BUC}(\mathbb{R})\). Standard properties of the heat semigroup imply that 
        \[ 
            J\in C\bigl([0,T];\operatorname{BUC}(\mathbb{R})\bigr) \cap C^{1,2}\bigl((0,T]\times\mathbb{R}\bigr), \qquad \partial_t J=k\partial_{xx}J, 
        \] 
        and 
        \[ 
            \lim_{t\downarrow0}\|J_t-J_0\|_\infty=0. 
        \] 
        Moreover, since \(V_0\geq0\), we have
        \[ 
            J_t(x) = (K_t*J_0)(x) \geq \inf_{x\in\R} J_0(x) \geq \cosh\rho>1. 
        \] 
        Consequently, we may define 
        \[ 
            V_t :=\frac{1}{\gamma} \left[\arccosh(J_t)-\rho\right] =\frac{1}{\gamma} \left[ \arccosh \left(K_t*\cosh(\rho+\gamma V_0)\right)-\rho \right]. 
        \] 
        The preceding lower bound immediately yields \(V_t\geq0\). Since \(J_t\) takes values in the compact interval 
        \[ 
            [\cosh\rho,\cosh(\rho+\gamma\|V_0\|_\infty)], 
        \] 
        and $\arccosh$ is smooth on $(1,\infty)$ and Lipschitz with constant $1/\sinh(\rho)$ on this interval, it follows that 
        \[ 
            V\in C\bigl([0,T];\operatorname{BUC}(\mathbb{R})\bigr) \cap C^{1,2}\bigl((0,T]\times\mathbb{R}\bigr) 
        \] 
        and 
        \[ 
            \|V_t-V_0\|_\infty \leq \frac{1}{\gamma\sinh\rho}\|J_t-J_0\|_\infty \longrightarrow0 \qquad\text{as }t\downarrow0. 
        \] 
        From the proof of Lemma~\ref{lem:genburg} it follows that $V_t$ satisfies 
        \[
            \partial_t V_t(x) = k \partial_{xx}V_t(x) + k\gamma \coth\bigl(\rho+\gamma V_t(x)\bigr) \bigl[\partial_xV_t(x)\bigr]^2, \qquad (t,x)\in(0,T]\times\mathbb{R}.
        \]
        For uniqueness, let \(\widetilde V_t\) be another solution with initial datum \(V_0\) and set \[ \widetilde J_t:=\cosh(\rho+\gamma\widetilde V_t). \] The same computation as in Lemma~\ref{lem:genburg} shows that \(\widetilde J_t\geq \cosh(\rho)\) is solution of 
        \[ 
            \partial_t\widetilde J_t=k\partial_{xx}\widetilde J_t,\qquad (t,x)\in(0,T]\times\mathbb{R},
        \]
        satisfying 
        \[
            \lim_{t\downarrow 0}||\widetilde J_t-J_0||_\infty =0. 
        \] 
        Uniqueness of solutions of the heat equation in this class of functions therefore gives 
        \[ 
            \widetilde J_t=K_t*J_0=J_t.
        \] 
        Hence, 
        \[ 
            \cosh(\rho+\gamma\widetilde V_t) =\cosh(\rho+\gamma V_t). 
        \] 
        The fact that \(\rho+\gamma V\geq\rho>0\) and injectivity of $\cosh$ on $[0,\infty)$ imply that \(\widetilde V_t=V_t\). This proves the result. \end{proof}
        Since in our context $V_t(z)=\int_z^\infty \nu_t(x)\dd x$ is the tail potential of some energy field, we expect it to be non-increasing in $z$ and satisfy 
        \begin{align*}
            \lim_{z\to-\infty} V_t(z) \leq E_\mathrm{D},\\
            \lim_{z\to\infty} V_t(z) = 0,
        \end{align*} 
        where $E_\mathrm{D}\geq 0$ is the bound on the total energy. The following proposition shows that the non-compact dynamic deterministic KPZ equation indeed preserves monotonicity of the initial condition and satisfies both these limits. Even more is true. If $\nu_0(z)$ is not zero a.e., then for each time $t>0$ we have $\nu_t(z)= -\partial_zV_t >0$. This infinite speed of propagation is inherited from the heat equation.
        \begin{proposition} \label{prop:propertiesgenburger}
            Let $T,k,\gamma,\rho >0$, $V_0\in C_b(\R)$ nonnegative and let
        \[ 
            V_t = \frac{1}{\gamma} \left[ \arccosh\left( K_t * \cosh\bigl(\rho+\gamma V_0\bigr)\right) -\rho \right] \geq 0.
        \] 
        Then the following statements hold.
        \begin{enumerate}[label=(\roman*)]
            \item If $V_0$ is nondecreasing, nonincreasing respectively, then $V_t$ is nondecreasing, nonincreasing respectively, for all $t\in[0,T]$.
            \item If $V_0$ is nonconstant and nondecreasing, nonincreasing respectively, then for each $t\in(0,T]$ and $z\in\R$ we have $\partial_z V_t(z)>0$, $\partial_z V_t(z)<0$ respectively.
            \item If $$\lim_{z\to \pm \infty} V_0(z) = L_{\pm}$$ for some constant $L_\pm\geq 0$, then $$\lim_{z\to \pm \infty} V_t(z) = L_{\pm}$$ for all $t\in [0,T]$.
        \end{enumerate}
        \end{proposition}
        \begin{remark}
            We only require $V_0 \in C_b(\R)$. However, in all three cases of the proposition we automatically have $V_0\in \BUCR$. Indeed, if a continuous function is either bounded and monotone or has finite limits at both $\pm \infty$, it is uniformly continuous as well. 
        \end{remark}
        \begin{proof} 
        Suppose first that $V_0$ is nondecreasing. Set \[ J_0(x)=\cosh\bigl(\rho+\gamma V_0(x)\bigr) \] and \[ J_t=K_t*J_0, \] so that
        \[
            V_t = \frac{\arccosh(J_t)-\rho}{\gamma}.
        \]
        Since $V_0\geq0$ and $\rho,\gamma>0$, the functions \[ v\longmapsto \cosh(\rho+\gamma v) \quad\text{and}\quad m\longmapsto \frac{\arccosh(m)-\rho}{\gamma} \] are strictly increasing on the relevant ranges. Therefore, $J_0$ is nondecreasing. For $x_1<x_2$, 
        \begin{align}
            J_t(x_2)-J_t(x_1) = \int_{\mathbb{R}}K_t(y) \bigl[J_0(x_2-y)-J_0(x_1-y)\bigr]\,dy \geq0. \label{eq:Mtx1-Mtx2}
        \end{align}
        Hence $J_t$, and therefore $V_t$, is nondecreasing. The nonincreasing case follows analogously. This proves~(i). 
        
        For (ii), fix $t>0$ and suppose $V_0$ is nonconstant and nondecreasing. Then $J_0$ is bounded, nonconstant and nondecreasing. Since $\partial_x\arccosh(x) > 0$, it is enough to show that 
        $$\partial_x J_t(x) = \partial_x \int_\R K_t(x-y) J_0(y) \dd y >0.$$ Let us justify that we can interchange the partial derivative and integral by showing that for fixed $x_0\in\R$, there exists an integrable function $f$ such that 
        \[
            |\partial_x K_t(x-y) J_0(y)| \leq f(y)
        \]
        for all $x\in (x_0-1,x_0+1)$. Since $$\partial_xK_t(x-y)=-\frac{x-y}{2kt}K_t(x-y),$$ we have for $|x-x_0|\leq 1$ that
        \[
            |\partial_xK_t(x-y)J_0(y)|\leq f(y),
        \]
        where 
        \[
            f(y)= (1+|x_0-y|)\|J_0\|_\infty(2kt\sqrt{4\pi kt})^{-1} \cdot \begin{cases}
                1 &\text{if } |y-x_0|<2, \\
                e^{-(y-x_0)^2/(16kt)}\quad &\text{else.}
            \end{cases}
        \]
        Here we used that $|x-y|\leq |x-x_0|+|x_0-y| \leq 1 + |x_0-y|$ and,  by the reversed triangle inequality, for $|x_0-y|\geq 2$ that
        $$|x-y|\geq |x_0-y| - |x_0-x| \geq |x_0-y| -1 \geq \frac12 |x_0-y|.$$ 
        Since $f\in L^1(\R)$, we can interchange derivative and integral, and after substituting $u=x-y$ and using that $K_t$ is even,
        \[
            \partial_xJ_t(x) = -\int_\R \frac{u}{2kt}K_t(u)J_0(x-u)\dd u
            = \int_0^\infty \frac{u}{2kt}K_t(u)\big[J_0(x+u)-J_0(x-u)\big]\dd u \geq 0.
        \]
        If $\partial_xJ_t(x)=0$, then $J_0(x+u)=J_0(x-u)$ for all $u>0$. By monotonicity, $J_0(x-u)\leq J_0(x)\leq J_0(x+u)$, so $J_0$ would be constant, a contradiction. Hence $\partial_xJ_t(x)>0$. The nonincreasing case is analogous.

    Finally, suppose that \[ \lim_{x\to\pm\infty}V_0(x)=L_\pm. \] Note that $J_0$ is bounded because $V_0$ is. Since we also have $||K_t||_1=1$, we can use the dominated convergence theorem to obtain
        \begin{align*}
            \lim_{x\to\pm \infty} J_t(x)= \lim_{x\to\pm \infty} \big[K_t * \cosh(\rho+\gamma V_0) \big](x) = \int_\R K_t(y)\cosh(\rho+\gamma L_\pm)\dd y = \cosh(\rho+\gamma L_\pm).
        \end{align*}
        Hence $$\lim_{x\to \pm \infty} V_t(x) = L_{\pm},$$ proving (iii). 
        \end{proof}
\section*{Acknowledgment}
We thank Tim Roozendaal for useful discussions. Part of the work on the hydrodynamic limit of ABEP was performed during his master project.

\end{document}